\documentclass[12pt]{article}

\usepackage{amssymb}

\begin{document}
\bigskip

\begin{center}
{\Large \bf
Conjugacy results for the Lie algebra $\mathfrak{sl}_{\bf 2}$ \\
\vspace{0.05cm}
over an algebra which is a UFD}
\end{center}
\bigskip

\ \bigskip

\begin{center}
{\large Stephen Berman\footnote{
Partially supported by an NSERC Grant}
}
\end{center}

\begin{center}
Department of Mathematics\\
University of Saskatchewan\\
Saskatoon, Saskatchewan, S7N 5E6\\
CANADA
\end{center}
\bigskip

\begin{center}
{\large Jun Morita\footnote{
Partially supported by a Monkasho Kakenhi (2003 -- 2004)}
}
\end{center}

\begin{center}
Institute of Mathematics\\
University of Tsukuba\\
Tsukuba, Ibaraki, 305-8571\\
JAPAN
\end{center}
\bigskip

\ \bigskip

\noindent
{\bf Abstract}.\quad
Let $F$ be a field of characteristic not $2$ and assume all algebras are over F. We  establish several conjugacy
theorems for the special linear Lie algebra
$\mathfrak{sl}_2$ over an $F$-algebra which is a UFD. 
We  find the structure of the full automorphism group 
as well as conditions for when two such algebras are isomorphic.
We  also study the structure of their derivation algebras.
Often we work with more general coordinates when this is feasible.
\bigskip

\ \bigskip

\bigskip

\ \bigskip


\newpage

\section{Introduction}

This paper is about $\mathfrak{sl}_2(R)$ where $R$ is a commutative associative $F$-algebra with identity which is an integral domain and 
where $F$ is a field of characteristic different from $2$. We often will assume that $R$ is also a unique factorization domain but we work in more generality when possible. We are interested in proving various conjugacy type results in this setting.
The conjugacy results we establish involve either certain types of elements or certain types of subalgebras of
$\mathfrak{sl}_2(R)$. Also, the conjugacy results either allow the conjugating automorphism to come from the full automorphism group or, in some cases, we are able to use a specific subgroup of the full automorphism group. Moreover, we wanted to approach this problem in a way that worked for all characteristics different from $2$ and one that was as elementary as possible. More often than not, our results are computational in nature, and they actually give the exact form of the conjugating automorphism rather than just existence results. Of course, to do this we are taking advantage of the fact that two-by-two matrices are fairly easy to compute with. However, the reader should appreciate that to pick the right matrix for the particular setting  one is working  in is usually not an easy matter.

   Recall from \cite{H},\cite{J},\cite{K} that all Cartan subalgebras of a finite dimensional simple Lie algebra over an algebraically closed field of characteristic zero are conjugate. In fact the conjugating automorphism may be choosen to be a
so called inner automorphism. This is one of the central results in the theory of the finite dimensional simple Lie algebras
and one that has several different proofs. Indeed the three proofs in the references \cite{H},\cite{J},\cite{K} are all quite different. One is algebraic, one relies on algebraic geometry, and one (just over the complex numbers) on analysis. Moreover, over fields which are not algebraically closed it is known that finite dimensional simple Lie algebras do not necessiarily have all of their Cartan subalgebras being conjugate, see \cite{K}. It is with this hindsite that there have been some recent investigations of subalgebras in some infinite dimensional Lie algebras which play the role of Cartan subalgebras and study their conjugacy properties. It is the papers \cite{P1},\cite{P2} which stirred our interest in this question. In these papers  
A. Pianzola has studied Lie algebras of the form $\mathfrak{g} \otimes R$ where $R$ is a commutative $F$-algebra for a field $F$ of characteristic zero. He has used an algebraic geometric tact in these studies and has proved some deep results on the conjugacy of the so called MAD's. A MAD is a maximal abelian ad-diagonalizable subalgebra of the Lie algebra. In this setting these MAD's play the role of Cartan subalgebras. In fact, it is easy to see that in a finite dimensional simple Lie algebra over an algebraically closed field of characteristic zero that a Cartan subalgebra and a MAD are one and the same object. From the beginning our feeling was that, at least in the case when $\mathfrak{g}= \mathfrak{sl}_2(F)$, there should be another proof of the conjugacy of MAD's  when $R$ was a UFD. Moreover, the proof should be algebraic, fairly simple, and work in any characteristic except two. We were also inspired by \cite{M} and found that, when we got going with our investigation, we were able to not only prove conjugacy results for MAD's in this setting but several other conjugacy type results as well. 

   Here is an outline of what can be found in the paper. In the section following this Introduction we begin with a key Lemma which plays a role throughout the paper. It deals with what can be thought of as the arithmetic of $s\ell_2$-triples 
$(X,H,Y)$ in $\mathfrak{sl}_2(R)$ satisfying
$$[H,X]=2X,\quad [H,Y]=-2Y,\quad [X,Y]=H.$$
In particular we show that the middle element, $H$, must satisfy $H^2=I$ where $I$ is the two-by-two identity matrix.
Our next result then shows that any $H$ satisfying $H^2=I$ is the middle element of some $s\ell_2$-triple when $R$ is a UFD.
The final result of this section shows that all such elements are conjugate in the UFD case.

   The next section is the heart of the paper. We begin by showing that for a UFD all MAD's are conjugate and then use this to determine the structure of the automorphism group of $\mathfrak{sl}_2(R)$. We show that
$Aut_F(\mathfrak{sl}_2(R))$ is isomorphic to the semi-direct product of $PGL_2(R)$ and $Aut_F(R)$. 
Next we show that, only assuming that $R$ is  an integral domain, that all $s\ell_2$-triples are conjugate and, moreover, if $R$ is a UFD then the conjugating automorphism may be taken from $PGL_2(R)$. We close this section by showing when $R$ is a UFD that all elements $X$ in the first place of an $s\ell_2$-triple are conjugate and then determine, when $F=F^2$, what all elements, which belong to an isomorphic copy of $\mathfrak{sl}_2(F)$ inside of $\mathfrak{sl}_2(R$), look like.

In the final section we show $\mathfrak{sl}_2(R)$ can be isomorphic to $\mathfrak{sl}_2(S)$ for $R$ an integral domain and $S$ a general commutative $F$-algebra only when $R$ and $S$ are isomorphic as $F$-algebras. We then go on to determine the derivations of 
$\mathfrak{sl}_2(R)$. The paper closes with several remarks. In these we either point out some other interesting facts about $\mathfrak{sl}_2(R)$ or compare our results with those in \cite{P1},\cite{P2}. The reader should remember that these two references assume that the characteristic of the underlying field is zero so we only present a comparison for characteristic zero fields.

\section{Basic results on $s\ell_2$-triples in $\mathfrak{sl}_2(R)$}

Throughout we let $F$ be a field of characteristic not $2$.
We always let $R$ denote a unital commutative associative 
$F$-algebra and often we will assume that $R$ is an integral domain or 
a UFD. When this is the case we just say {R} is an integral domain or a UFD so 
the reader should always understand that such an $R$ is also an $F$- algebra. We 
identify $F$ as a subset of $R$ and let $R^\times$ be the unit group of $R$. When $R$ is 
an integral domain we let $K$ denote it's quotient field so that $K$ is an extension field of $F$. 

We put
$$\mathfrak{sl}_2(R) = \left\{
\ \left. \left(
\begin{array}{cc}
a & b\\ c & -a
\end{array}
\right) \ \right| \ 
a,b,c,d \in R\ 
\right\} \ .$$
Sometimes $\mathfrak{sl}_2(R)$ is regarded as
a subspace of the full matrix algebra
$M_2(R)$ or as an $F$-subspace of $\mathfrak{sl}_2(K)$.
Of course we have that $\mathfrak{sl}_2(F)$ is a subalgebra of $\mathfrak{sl}_2(R)$.
Put
$${\mathfrak h} = F \left(
\begin{array}{cc}
1 & 0\\ 0 & -1
\end{array}
\right) \ .$$
This is a one dimensional subalgebra of $\mathfrak{sl}_2(R)$.

Recall that a triple $(e,h,f)$ of nonzero elements of a Lie algebra
is called an $s\ell_2$-triple
if
$$[h,e]=2e,\quad [h,f]=-2f,\quad [e,f]=h.$$
Then the $F$-span of the three elements $e,h,f$ is isomorphic 
to the Lie algebra $\mathfrak{sl}_2(F)$. Notice that if $(e,f,h)$ is 
an $s\ell_2$-triple then so is $(f,-h,e)$. Moreover, if $(X,H,Y)$
is an $s\ell_2$-triple in $\mathfrak{sl}_2(R)$ then so is
$(Y^t,H^t,X^t)$ where for $Z$ an element of $\mathfrak{sl}_2(R)$ 
we are letting $Z^t$ denote the transpose of $Z$. Also, we let $I$ denote the 
identity matrix in $M_2(R)$.
   We want to study $s\ell_2$-triples in $\mathfrak{sl}_2(R)$ and to do so
we take three elements:
$$
X =
\left(
\begin{array}{cc}
X_1 & X_2\\ X_3 & -X_1
\end{array}
\right) \ ,
\qquad
H =
\left(
\begin{array}{cc}
H_1 & H_2\\ H_3 & -H_1
\end{array}
\right) \ ,
\qquad
Y =
\left(
\begin{array}{cc}
Y_1 & Y_2\\ Y_3 & -Y_1
\end{array}
\right)
$$
in $\mathfrak{sl}_2(R)$.

\bigskip

\noindent
{\bf Lemma 1}. Suppose that $R$ is an integral domain and
$(X,H,Y)$ is an $s\ell_2$-triple with $X,H,Y \in \mathfrak{sl}_2(R)$.
Then:\par
\noindent
(1)\quad
$H_1^2 + H_2H_3 = 1$,\par
\noindent
(2)\quad
$X_1^2 + X_2X_3 = 0$,\par
\noindent
(3)\quad
$Y_1^2 + Y_2Y_3 = 0$,\par
\noindent
(4)\quad 
$H^2 =I
$.
\bigskip

{\it Proof}\quad  Since $(X,H,Y)$ is an $s\ell_2$-triple, we have
$$[H,X]=2X,\quad [H,Y]=-2Y,\quad [X,Y]=H,$$
which implies
$$\begin{array}{lllll}
{\rm (A)} & \left\{ 
\begin{array}{llll}
{\rm (A1)} & \quad H_2X_3 - H_3X_2 & = & 2X_1,\\
{\rm (A2)} & \quad H_1X_2 - H_2X_1 & = & X_2,\\
{\rm (A3)} & \quad H_3X_1 - H_1X_3 & = & X_3,
\end{array}
\right.\\
&\\
{\rm (B)} & \left\{ 
\begin{array}{llll}
{\rm (B1)} & \quad H_2Y_3 - H_3Y_2 & = & -2Y_1,\\
{\rm (B2)} & \quad H_1Y_2 - H_2Y_1 & = & -Y_2,\\
{\rm (B3)} & \quad H_3Y_1 - H_1Y_3 & = & -Y_3,
\end{array}
\right.\\
&\\
{\rm (C)} & \left\{ 
\begin{array}{llll}
{\rm (C1)} & \quad X_2Y_3 - X_3Y_2 & = & H_1,\\
{\rm (C2)} & \quad 2(X_1Y_2 - X_2Y_1) & = & H_2,\\
{\rm (C3)} & \quad 2(X_3Y_1 - X_1Y_3) & = & H_3.
\end{array}
\right.
\end{array}$$
We can rewrite (A) as
$$
\left\{
\begin{array}{lrcrcrcl}
{\rm (A1)}' & \quad -2X_1 & - & H_3X_2 & + & H_2X_3 & = & 0,\\
{\rm (A2)}' & \quad -H_2X_1 & + & (H_1-1)X_2 &&& = & 0,\\
{\rm (A3)}' & \quad H_3X_1 &&& - & (H_1+1)X_3 & = & 0.
\end{array}
\right.
$$
Then, from ${\rm (A1)}' \times (-H_2) + {\rm (A2)}' \times 2$, we obtain
$$(H_2H_3+2H_1-2)X_2 - H_2^2X_3 = 0$$
and, from ${\rm (A1)}' \times H_3 + {\rm (A3)}' \times 2$, we obtain
$$-H_3^2X_2 + (H_2H_3-2H_1-2)X_3 = 0.$$
Therefore, we can write
$$\left(
\begin{array}{cc}
H_2H_3+2H_1-2 & -H_2^2\\
-H_3^2 & H_2H_3-2H_1-2
\end{array}
\right)
\ 
\left(
\begin{array}{c}
X_2\\ X_3 \end{array}
\right)
=
\left(
\begin{array}{c}
0\\ 0 \end{array}
\right) \ .
$$
Since $X$ is a nonzero element, 
we see
$$\left(
\begin{array}{l} X_2\\ X_3 \end{array}
\right)$$
is a nonzero vector by ${\rm (A1)}$.
Hence, we have
$$\begin{array}{lll}
{\rm det}\ 
\left(
\begin{array}{cc}
H_2H_3+2H_1-2 & -H_2^2\\
-H_3^2 & H_2H_3-2H_1-2
\end{array}
\right)
& = & (H_2H_3 - 2)^2 - 4H_1^2 - H_2^2H_3^2\\
& = & -4H_2H_3 + 4 - 4H_1^2\\
& = & 0
\end{array}
$$
and $H_1^2 + H_2H_3 = 1$, which is (1). Then,
$$H^2 =
\left(
\begin{array}{cc}
H_1^2+H_2H_3 & 0\\
0 & H_1^2+H_2H_3
\end{array}
\right)
=
\left(
\begin{array}{cc}
1 & 0\\
0 & 1
\end{array}
\right) \ ,
$$
which shows (4).
Computing ${\rm (A2)} \times X_3 + {\rm (A3)} \times X_2$, we obtain
$$-H_2X_1X_3 + H_3X_1X_2 = 2X_2X_3.$$
Hence we get, by ${\rm (A1)}$,
$$-2X_1^2 = 2X_2X_3$$
and $X_1^2 + X_2X_3 = 0$, which is (2). Similarly, using ${\rm (B)}$, we can
obtain $Y_1^2 + Y_2Y_3 = 0$ and this is (3).
Q.E.D.
\bigskip

   In our next result we characterize elements $H$ satisfying $(4)$ in
Proposition 1. For this we will need to assume that R is a UFD. We keep the 
same notation as above.
\bigskip

\noindent
{\bf Proposition 2}. Suppose that $R$ is a UFD.
If $H$ is an element of $\mathfrak{sl}_2(R)$ with
$$
H^2 =I
$$
then there exist $X,\ Y \in \mathfrak{sl}_2(R)$ such that
$(X,H,Y)$ is an $s\ell_2$-triple.
\bigskip

{\it Proof}\quad  We will give an explicit
construction and will show where we are using the UFD property of $R$.
First we suppose $H_2 = 0$, which implies $H_1 = \pm 1$ by $(1.)$ of Proposition 1.
If $H_1 = 1$, we put
$$X = \left(
\begin{array}{cc}
\displaystyle{ - \frac{H_3}{2} } & 1\\
&\\
\displaystyle{ - \frac{H_3^2}{4} } & \displaystyle{ \frac{H_3}{2} }
\end{array}
\right)
\qquad \mbox{and}\qquad
Y = \left(
\begin{array}{cc}
0 & 0\\
&\\
1 & 0
\end{array}
\right) \ .
$$
Then, $(X,H,Y)$ is an $s\ell_2$-triple.
If $H_1 = -1$, let
$$X = \left(
\begin{array}{cc}
0 & 0\\
&\\
-1 & 0
\end{array}
\right)
\qquad \mbox{and}\qquad
Y = \left(
\begin{array}{cc}
\displaystyle{ - \frac{H_3}{2} } & -1\\
&\\
\displaystyle{ \frac{H_3^2}{4} } & \displaystyle{ \frac{H_3}{2} }
\end{array}
\right) \ .
$$
Then, we also see that $(X,H,Y)$ is an $s\ell_2$-triple.
Indeed, this follows from the first computation above by replacing $H_3$
by $-H_3$ and then replacing the $s\ell_2$-triple we obtain, say it is (e,h,f),
by the $s\ell_2$-triple (-f,-h, -e).
\bigskip

\noindent
Next we suppose $H_3 = 0$, which implies $H_1 = \pm 1$.
In the case when $H_1 = 1$, we put
$$X = \left(
\begin{array}{cc}
0 & 1\\
&\\
0 & 0
\end{array}
\right)
\qquad \mbox{and}\qquad
Y = \left(
\begin{array}{cc}
- \displaystyle{ \frac{H_2}{2} } & - \displaystyle{ \frac{H_2^2}{4} }\\
&\\
1 & \displaystyle{ \frac{H_2}{2} }
\end{array}
\right) \ .
$$
Then, $(X,H,Y)$ is an $s\ell_2$-triple. Indeed, this is 
just the transpose of the first relation above after we replace $H_3$ 
by $H_2$. 
In the case  $H_1 = -1$, we put
$$X = \left(
\begin{array}{cc}
- \displaystyle{ \frac{H_2}{2} } & \displaystyle{ \frac{H_2^2}{4} }\\
&\\
-1 & \displaystyle{ \frac{H_2}{2} }
\end{array}
\right)
\qquad \mbox{and}\qquad
Y = \left(
\begin{array}{cc}
0 & -1\\
&\\
0 & 0
\end{array}
\right) \ .
$$
Then, we also see that $(X,H,Y)$ is an $s\ell_2$-triple. Indeed, this follows from the 
third case above just as the second followed from the first.
\bigskip

\noindent
Finally we suppose that $H_2 \not= 0$ and $H_3 \not= 0$.
This means $H_1 \not= \pm 1$
since
$$(H_1-1)(H_1+1) = - H_2H_3.$$

 It is here that we use our hypothesis that $R$ is a UFD.
We write $gcd(x,y)$ for the greatest common divisor of two elements of $R$.
Notice that $gcd(H_1-1,H_1+1)=1$  since $(H_1+1) - (H_1-1) = 2 \in R^\times$.
Let $a=gcd(H_1-1,H_2)$ and $a'=gcd(H_1+1, H_2)$ so that $gcd(a,a')=1$. Also let 
$b=gcd(H_1-1,H_3)$ and $b'=gcd(H_1+1,H_3)$ which means we have 
$$gcd(b,b')=gcd(a,b')=gcd(a',b)=1.$$
Thus $aa'$ divides $H_2$ and $bb'$ divides $H_3$, hence the product
$aa'bb'$ divides $H_2H_3=-(H_1-1)(H_1+1)$.
Clearly $gcd(a,H_1+1)=1$ and $gcd(b,H_1+1)=1$ so we obtain that  the product $ab$ divides
$H_1-1$. Similiarly $a'b'$ divides $H_1+1$.

   Now if $ab$ and $H_1-1$ are not associates in $R$ there exists a prime $p$ of $R$ 
such that $abp$ divides $H_1-1$. Then $abp$ divides $H_2H_3$ so either $ap$ divides $H_2$ or 
$bp$ divides $H_3$. In either case we get a contradiction to the definition of $a$ or $b$.  

Thus
we can set
$${\rm (H)}\quad 
\left\{
\quad 
\begin{array}{lll}
H_1-1 & = & uab,\\
H_1+1 & = & u'a'b',\\
H_2 & = & vaa',\\
H_3 & = & wbb',
\end{array}
\right.$$
where $u,u',v,w \in R^\times$ with $uu' = -vw$,
and $a,a',b,b'$ are as above. 

In this situation we put
$$X' =
\left(
\begin{array}{cc}
(H_1-1)H_2 & H_2^2\\
-(H_1-1)^2 & -(H_1-1)H_2
\end{array}
\right)
$$
and
$$
Y' =
\left(
\begin{array}{cc}
-(H_1+1)H_2 & -H_2^2\\
(H_1+1)^2 & (H_1+1)H_2
\end{array}
\right)
$$
as nonzero elements of $\mathfrak{sl}_2(R) \subset \mathfrak{sl}_2(K)$.
(Recall that $K$ is the quotient field of $R$).
Then, we easily find the following relations:
$$[ H , X' ] = 2X' \qquad \mbox{and}\qquad [ H , Y' ] = -2Y',$$
which gives
$$\mathfrak{sl}_2(K) = KX' \oplus KH \oplus KY'$$
as the eigenspace decomposition of ${\rm ad}\ H$ on
$\mathfrak{sl}_2(K)$. Here $KX'$, $KH$ and $KY'$ are the eigenspaces
corresponding to $2,0,-2$ respectively. In particular we now choose
$X \in KX'$ and $Y \in KY'$ as follows.
Let
$$X = \displaystyle{\frac{1}{2va^2}} X'
\qquad \mbox{and}\qquad
Y = \displaystyle{\frac{1}{2va'^2}} Y'\ .$$
It is clear that these are nonzero elements of $\mathfrak{sl}_2(R)$.
Then, by a direct calculation, we have $[X,Y] = H$. Therefore,
$(X,H,Y)$ is our desired  $s\ell_2$-triple. Q.E.D.
\bigskip

Recall that $GL_2(R)$ is the group of two by two matrices with entries 
in $R$ which are invertible with their inverses also in $GL_2(R)$. Thus 
this is the group of matrices in $M_2(R)$ which have their determinants 
in $R^\times$. If $P \in GL_2(R)$ then the map $ X \rightarrow P^{-1}XP$
for all $X \in \mathfrak{sl}_2(R)$ is an automorphism of $\mathfrak{sl}_2(R)$.

Our first two results have shown that elements $H \in \mathfrak{sl}_2(R)$ which
satisfy $H^2=I$   are exactly the elements which are the middle elements 
in an $s\ell _2$-triple. We next show all such elements are conjugate by elements from $GL_2(R)$.
Our method of proof is similar to that used in Proposition 2 and, in fact, we will use some of the notation developed there.
\bigskip

\noindent
{\bf Theorem 3}. Suppose that $R$ is a unique factorization domain.
If $H$ is an element of $\mathfrak{sl}_2(R)$ with
$$
H^2 = I
$$
then there is an element $P \in GL_2(R)$ such that
$$P^{-1}HP =
\left(
\begin{array}{cc}
1 & 0\\ 0 & -1
\end{array}
\right)\ .
$$
\bigskip
{\it Proof}\quad  We take
$$H = \left( 
\begin{array}{cc}
H_1 & H_2\\ H_3 & -H_1
\end{array} \right)
\in \mathfrak{sl}_2(R).$$
Suppose 
$$H^2 = \left(
\begin{array}{cc}
H_1^2+H_2H_3 & 0\\
0 & H_1^2+H_2H_3
\end{array}
\right)
=  I
$$
so that
$H_1^2 + H_2H_3 = 1$.
Then, as in the previous result,
we have the following three cases:
\bigskip

\noindent
(Case 1)\quad $H_1 = \pm 1$ and $H_2 = 0$.\par
\noindent
(Case 2)\quad $H_1 = \pm 1$ and $H_3 = 0$.\par
\noindent
(Case 3)\quad $H_1 \not= \pm 1$, $H_2 \not= 0$ and $H_3 \not= 0$.
\bigskip

\noindent
\underline{(Case 1)}:
\bigskip

\noindent
If $H_1 = 1$, then we have
$$
\left(
\begin{array}{cc}
1 & 0\\ - \displaystyle{ \frac{H_3}{2} } & 1
\end{array}
\right)
\ 
\left(
\begin{array}{cc}
1 & 0\\ H_3 & - 1
\end{array}
\right)
\ 
\left(
\begin{array}{cc}
1 & 0\\ \displaystyle{ \frac{H_3}{2} } & 1
\end{array}
\right)
=
\left(
\begin{array}{cc}
1 & 0\\ 0 & - 1
\end{array}
\right) \ .$$
If $H_1 = - 1$, then we have
$$
\left(
\begin{array}{cc}
\displaystyle{ \frac{H_3}{2} } & 1\\ - 1 & 0
\end{array}
\right)
\ 
\left(
\begin{array}{cc}
- 1 & 0\\ H_3 & 1
\end{array}
\right)
\ 
\left(
\begin{array}{cc}
0 & - 1\\ 1 & \displaystyle{ \frac{H_3}{2} }
\end{array}
\right)
=
\left(
\begin{array}{cc}
1 & 0\\ 0 & - 1
\end{array}
\right) \ .$$
\bigskip

\noindent
\underline{(Case 2)}:
\bigskip

\noindent
If $H_1 = 1$, then we have
$$
\left(
\begin{array}{cc}
1 & \displaystyle{ \frac{H_2}{2} }\\ 0 & 1
\end{array}
\right)
\ 
\left(
\begin{array}{cc}
1 & H_2\\ 0 & - 1
\end{array}
\right)
\ 
\left(
\begin{array}{cc}
1 & - \displaystyle{ \frac{H_2}{2} }\\ 0 & 1
\end{array}
\right)
=
\left(
\begin{array}{cc}
1 & 0\\ 0 & - 1
\end{array}
\right) \ .$$
If $H_1 = - 1$, then we have
$$
\left(
\begin{array}{cc}
0 & 1\\ -1 & \displaystyle{ \frac{H_2}{2} }
\end{array}
\right)
\ 
\left(
\begin{array}{cc}
- 1 & H_2\\ 0 & 1
\end{array}
\right)
\ 
\left(
\begin{array}{cc}
\displaystyle{ \frac{H_2}{2} } & -1\\ 1 & 0
\end{array}
\right)
=
\left(
\begin{array}{cc}
1 & 0\\ 0 & - 1
\end{array}
\right) \ .$$
\bigskip

\noindent
\underline{(Case 3)}:
\bigskip

\noindent
Put
$$V = \left(
\begin{array}{c}
H_1+1\\ H_3
\end{array}
\right)
\qquad \mbox{and}\qquad
W = \left(
\begin{array}{c}
H_1-1\\ H_3
\end{array}
\right) \ .
$$
Then, $HV = V$ and $HW = -W$.
We write, as in the previous result,
$$\begin{array}{lll}
H_1 - 1 & = & uab,\\
H_1 + 1 & = & u'a'b',\\
H_3 & = & wbb'.
\end{array}$$
Put
$$\begin{array}{l}
V' = \displaystyle{ \frac{V}{b'} } =
\left( \begin{array}{cc} u'a' \\ wb \end{array} \right)\ ,\\
\\
W' = \displaystyle{ \frac{W}{b} } =
\left( \begin{array}{cc} ua \\ wb' \end{array} \right)\ ,\\
\\
P = \left(
\begin{array}{cc}
u'a' & ua\\
wb & wb'
\end{array} \right) \ .
\end{array}$$
Then
$$H P = P\ 
\left(
\begin{array}{cc}
1 & 0\\ 0 & -1
\end{array} \right)
$$
and ${\rm det}\ P = w(u'a'b' - uab) = w((H_1+1) - (H_1-1)) = 2w 
\in R^\times$.
In particular we see $P \in GL_2(R)$. Q.E.D.
\bigskip

\bigskip

\section{ Conjugacy results and 
automorphisms of $\mathfrak{sl}_2(R)$}

In the study of the finite dimensional simple Lie algebras over algebraically closed fields of characteristic zero one of the most important results is that all Cartan subalgebras are conjugate,
that is, there is an automorphism (one even of a certain specified form) taking one to the other. Our Theorem 3 in the previous section is also a conjugacy type result in so far as it relates all elements 
$H$ satisfying $H^2=I$ with one another. In this section we investigate various conjugacy results. Along the way we also obtain the structure of the automorphism group, $Aut_F(s\ell_2(R))$.

We will work with abelian subalgebras ${\mathfrak a}$. Recall that a subalgebra  ${\mathfrak a}$ of $\mathfrak{sl}_2(R)$ is
called ad-diagonalizable if ${\rm ad}\ Z$ is diagonalizable over $F$ on
$\mathfrak{sl}_2(R)$ for all $Z \in {\mathfrak a}$, that is,
for each $Z \in {\mathfrak a}$, the entire Lie algebra
$\mathfrak{sl}_2(R)$ has an $F$-basis, dependent of $Z$,
consisting of eigenvectors
of ${\rm ad}\ Z$ with eigenvalues in $F$.
More precisely, if we set
$$\mathfrak{sl}_2(R)_\alpha = \{ \ Z \in \mathfrak{sl}_2(R)\ \mid\ 
[A,Z] = \alpha(A) Z\ \mbox{for all}\ A \in {\mathfrak a}\ \},$$
for $\alpha \in {\mathfrak a}^* = Hom_F({\mathfrak a},F)$,
then that ${\mathfrak a}$ is ad-diagonalizable  means
$$\mathfrak{sl}_2(R) = 
\bigoplus_{\alpha \in {\mathfrak a}^*}\ \mathfrak{sl}_2(R)_\alpha.$$
Let ${\mathcal F}$ be the family of all
abelian ad-diagonalizable subalgebras of $\mathfrak{sl}_2(R)$.
An element ${\mathfrak a}$ of ${\mathcal F}$ is called a MAD (=
Maximal Abelian ad-Diagonalizable subalgebra) if
${\mathfrak a}$ is maximal in ${\mathcal F}$ with respect to inclusion.
The reader should understand that a MAD is one possible replacement for the notion of a Cartan 
subalgebra in our setting, see \cite{P1} \cite{P2}.
The next result shows they are all conjugate to each other.
\bigskip

\noindent
{\bf Theorem 4} Assume that $R$ is a unique factorization domain.
If ${\mathfrak h}'$ is a MAD of $\mathfrak{sl}_2(R)$, then
${\rm dim}_F\ {\mathfrak h}' = 1$ and
there is $P \in GL_2(R)$ such that
$P^{-1} {\mathfrak h}' P = {\mathfrak h}$.
\bigskip

\noindent
{\it Proof}\quad
Let
$$Z = \left(
\begin{array}{cc}
Z_1 & Z_2\\ Z_3 & - Z_1
\end{array}
\right)
$$
be a nonzero element of $\mathfrak{sl}_2(R)$. Then one of
the following three cases must hold:\par
\noindent
(Case 1) $Z_1^2 + Z_2Z_3 = 0$,\par
\noindent
(Case 2) $Z_1^2 + Z_2Z_3 \not= 0,\ Z_2 = 0$,\par
\noindent
(Case 3) $Z_1^2 + Z_2Z_3 \not= 0,\ Z_2 \not= 0$.\par
\noindent
We first develope some information on ${\rm ad}\ Z$ in each case.
\bigskip

\noindent
\underline{(Case 1)}: In this case, $Z^2 = O$. Therefore,
$({\rm ad}\ Z)^3 = 0$ and ${\rm ad}\ Z$
is nilpotent.
\bigskip

\noindent
\underline{(Case 2)}: In this case $Z_1\not=0$ and we put
$$P' = \left(
\begin{array}{cc}
1 & 0\\ \displaystyle{ \frac{Z_3}{2Z_1} } & 1
\end{array}
\right) \ .
$$
Then we have $P' \in GL_2(K)$ and
$$P'^{-1} Z P' =
\left(
\begin{array}{cc}
Z_1 & 0\\ 0 & - Z_1
\end{array}
\right) \ .
$$
In particular it follows that ${\rm ad}\ Z$ is diagonalizable on $\mathfrak{sl}_2(K)$.
\bigskip

\noindent
\underline{(Case 3)}:
Let $\bar K$ be the algebraic closure of $K$ and  choose 
a nonzero element $\lambda \in \bar K$
satisfying $\lambda^2 = Z_1^2 + Z_2Z_3$. Put
$$P'' = \left(
\begin{array}{cc}
1 & 1\\
&\\
- \displaystyle{ \frac{Z_1 - \lambda}{Z_2} } &
- \displaystyle{ \frac{Z_1 + \lambda}{Z_2} }
\end{array}
\right) \ .
$$
then we have
$$P''^{-1} Z P'' =
\left(
\begin{array}{cc}
\lambda & 0\\ 0 & - \lambda
\end{array}
\right) \ .
$$
Hence ${\rm ad}\ Z$ is diagonalizable on
$\mathfrak{sl}_2(\bar K)$.
\bigskip

\noindent
Now let ${\mathfrak h}'$ be a MAD of $\mathfrak{sl}_2(R)$.
We take a nonzero element $H' \in {\mathfrak h}'$
and note that it is an easy exercise to see that $\mathfrak{sl}_2(R)$ has no center.
Since ${\rm ad}\ H'$ is diagonalizable ${\rm ad}\ H'$ is
never nilpotent, otherwise ${\rm ad}\ H' = 0$, which means $H' = O$
and gives a contradiction. Therefore, by the three cases above, there is a nonzero element
$\lambda \in \bar K$
and a matrix $Q \in GL_2(\bar K)$ such that
$$Q^{-1} H' Q = \left(
\begin{array}{cc}
\lambda & 0\\ 0 & - \lambda
\end{array}
\right) \ .
$$
Hence, ${\rm ad}\ H'$ has eigenvalues $0,2\lambda,-2\lambda$
on $\mathfrak{sl}_2(\bar K)$, where $\lambda \in \bar K$. 
On the other hand, since $H'$ is in a MAD, we have
$$\mathfrak{sl}_2(R) = \oplus_{i \in I} FV_i$$
with $V_i \in \mathfrak{sl}_2(R)$ satisfying
$$[ H' , V_i ] = \mu_i V_i$$
for some $\mu_i \in F$, where $I$ is a suitable index set. 
In particular we see
$\mu_i = 0,\pm 2 \lambda$ for each $i \in I$. Not all $\mu_i $ can be zero so this   shows
$\lambda \in F$. We set
$$H'' = \frac{1}{\lambda}H'.$$
Then we have $H'' \in \mathfrak{sl}_2(R)$ and
$$H''^2 = \left(
\begin{array}{cc}
1 & 0\\ 0 & 1
\end{array}
\right) \ .
$$
Thus, by Theorem 3, we can find $P \in GL_2(R)$
such that
$$P^{-1} H'' P = \left(
\begin{array}{cc}
1 & 0\\ 0 & -1
\end{array}
\right) \ ,
$$
which implies
$$P^{-1} H' P = \left(
\begin{array}{cc}
\lambda & 0\\ 0 & - \lambda
\end{array}
\right) \ .
$$
Put ${\mathfrak h}'' = P^{-1} {\mathfrak h}' P$.
This is another MAD of $\mathfrak{sl}_2(R)$
containing ${\mathfrak h}$. Since ${\mathfrak h}''$ is abelian,
we obtain
$${\mathfrak h}'' \subset C_{\mathfrak{sl}_2(R)}({\mathfrak h})
= R
\left(
\begin{array}{cc}
1 & 0\\ 0 & -1
\end{array}
\right) \ .
$$
If we take a nonzero element
$$Z = \left(
\begin{array}{cc}
Z_1 & 0\\ 0 & - Z_1
\end{array}
\right) \in {\mathfrak h}'',
$$
then ${\rm ad}\ Z$ acts as the scalar $2Z_1$ on
$$\left(
\begin{array}{cc}
0 & R\\ 0 & 0
\end{array}
\right) \ .
$$
Indeed, we have
$$[ Z ,
\left(
\begin{array}{cc}
0 & 1\\ 0 & 0
\end{array}
\right) ]
=
\left(
\begin{array}{cc}
0 & 2Z_1\\ 0 & 0
\end{array}
\right) \ .
$$
Therefore, we obtain $2Z_1 \in F$, since $2Z_1$ is
an eigenvalue of ${\rm ad}\ Z$ on $\mathfrak{sl}_2(R)$.
Hence, $Z_1 \in F$, which shows ${\mathfrak h}'' \subset
{\mathfrak h}$. This implies ${\mathfrak h}'' = {\mathfrak h}$.
Thus, we have ${\rm dim}_F {\mathfrak h}' = 1$ and
$P^{-1} {\mathfrak h}' P = {\mathfrak h}$.
Q.E.D.
\bigskip

   We next want to use the above conjugacy result to determine the 
automorphism group of $\mathfrak{sl}_2(R)$ when $R$ is a UFD. For this first note 
that if $\rho$ is in $Aut_F(R)$ then $\rho$ acts on $\mathfrak{sl}_2(R)$ by acting on each 
matrix entry. In fact, it acts on $M_2(R)$ and $GL_2(R)$ as well.  Moreover, it is clear that this action induces an injective homomorphism from
 $Aut_F(R)$ to $Aut_F(\mathfrak{sl}_2(R))$. In this way we consider $Aut_F(R)$ as a subgroup 
of $Aut_F(\mathfrak{sl}_2(R))$.

Recall that if if $P \in GL_2(R)$ then the map, $ X \rightarrow PXP^{-1}$
for all $X \in \mathfrak{sl}_2(R)$, is an automorphism of $\mathfrak{sl}_2(R)$. We denote this automorphism by $\tau_P$. Moreover, if two elements, say $P$ and $Q$, 
of $GL_2(R)$ differ by an element $uI$ where $u \in R^\times$ then they induce the same automorphism.
That is, if $Q=uI P$ then the above conjugation action is the same for $P$ and $Q$ so $\tau_P =\tau_Q$.
 Now it is easy to see that the center, $ZGL_2(R)$ of 
$GL_2(R)$ is just $\{uI | u \in R^\times \}$. We let, as usual, the corresponding factor group be denoted by $PGL_2(R)$.
That is,
$$ PGL_2(R) = GL_2(R)/ZGL_2(R). $$
\noindent 
It is clear that the map $P \rightarrow \tau_P$ of $GL_2(R)$ to $Aut_F(\mathfrak{sl}_2(R)$ is a group 
homomorphism whose kernel is just $ZGL_2(R).$ Hence, in this way, we will consider $PGL_2(R)$ to be a subgroup of 
$Aut_F(\mathfrak{sl}_2(R)$. Finally note that this subgroup is normalized by the action of $Aut_F(R)$ since if
$\rho \in Aut_F(R)$ and $P \in GL_2(R)$ then we have 
$$\rho \tau_P \rho^{-1}= \tau_{\rho P}.$$ 
\bigskip

\noindent
{\bf Theorem 5}. For
$R$  a unique factorization domain we have
$$Aut_F(\mathfrak{sl}_2(R)) = PGL_2(R) \rtimes Aut_F(R).$$
\bigskip

\noindent
{\it Proof}\quad
Let $\sigma \in
Aut_F(\mathfrak{sl}_2(R))$ and  put
$$X = \sigma(
\left( \begin{array}{cc} 0 & 1\\ 0 & 0 \end{array} \right)
),\quad
H = \sigma(
\left( \begin{array}{cc} 1 & 0\\ 0 & -1 \end{array} \right)
),\quad
Y = \sigma(
\left( \begin{array}{cc} 0 & 0\\ 1 & 0 \end{array} \right)
).$$
Since $(X,H,Y)$ is an $sl_2$-triple, by Lemma 1, we have
$H^2 = I.$
Then, by using Theorem 3 we see 
there is $\tau' \in PGL_2(R)$ such that
$$\tau'(H) =
\left( \begin{array}{cc} 1 & 0\\ 0 & -1 \end{array} \right) \ .$$
Since
$$\begin{array}{lll}
[
\left( \begin{array}{cc} 1 & 0\\ 0 & -1 \end{array} \right)
,\tau' (X) ] & = & 2 \tau'(X),\\
&&\\
\rm {and } \quad [ \left( \begin{array}{cc} 1 & 0\\ 0 & -1 \end{array} \right),\tau' (Y) ] & = & -2 \tau'(Y),\\
&&\\ \rm {and } \quad
\left[ \tau' (X),\tau' (Y) \right] & = & 
\left(
\begin{array}{cc}
1 & 0\\ 0 & -1
\end{array} \right) \ ,
\end{array}
$$
we obtain
$$\begin{array}{lll}
\tau'(X) & = &
\left( \begin{array}{cc} 0 & u\\ 0 & 0 \end{array} \right) \\
&&\\
\tau'(Y) & = &
\left( \begin{array}{cc} 0 & 0\\ u^{-1} & 0 \end{array} \right)
\end{array}
$$
for some $u \in R^\times$. We choose $\tau'' \in PGL_2(R)$
corresponding to conjugation by
$$\left( \begin{array}{cc} u^{-1} & 0\\ 0 & 1 \end{array} \right) \ .
$$
Put $\tau = \tau'' \tau' \sigma$. Then, we have
$$\begin{array}{lll}
\tau(X)
& = & \left( \begin{array}{cc} 0 & 1\\ 0 & 0 \end{array} \right) \ ,\\
&&\\
\tau(H)
& = & \left( \begin{array}{cc} 1 & 0\\ 0 & -1 \end{array} \right) \ ,\\
&&\\
\tau(Y)
& = & \left( \begin{array}{cc} 0 & 0\\ 1 & 0 \end{array} \right) \ .
\end{array}$$
For each $r \in R$, we find
$$\tau(
\left( \begin{array}{cc} r & 0\\ 0 & -r \end{array} \right) )
= \left( \begin{array}{cc} r' & 0\\ 0 & -r' \end{array} \right)
$$
for some $r' \in R$, since
$$\begin{array}{lll}
[ \left( \begin{array}{cc} 1 & 0\\ 0 & -1 \end{array} \right) ,
\tau( 
\left( \begin{array}{cc} r & 0\\ 0 & -r \end{array} \right) ) ]
& = &
[ \tau( \left( \begin{array}{cc} 1 & 0\\ 0 & -1 \end{array} \right) ) ,
\tau(
\left( \begin{array}{cc} r & 0\\ 0 & -r \end{array} \right) ) ]
\\
&&\\
& = &
\tau([ \left( \begin{array}{cc} 1 & 0\\ 0 & -1 \end{array} \right) ,
\left( \begin{array}{cc} r & 0\\ 0 & -r \end{array} \right) ]  )
\\
&&\\
& = &
\left( \begin{array}{cc} 0 & 0\\ 0 & 0 \end{array} \right)
\end{array}$$
and
$$\tau(
\left( \begin{array}{cc} r & 0\\ 0 & -r \end{array} \right)
) \in C_{\mathfrak{sl}_2(R)}(
\left( \begin{array}{cc} 1 & 0\\ 0 & -1 \end{array} \right)
) = R\ 
\left( \begin{array}{cc} 1 & 0\\ 0 & -1 \end{array} \right) \ .
$$
Using this, we can define a map $\eta_\tau : R \longrightarrow R$
by $\eta_\tau(r) = r'$. 
Since
$$\begin{array}{lll}
\tau(
\left( \begin{array}{cc} r_1+r_2 & 0\\
0 & -(r_1+r_2) \end{array} \right) )
& = & \tau(
\left( \begin{array}{cc} r_1 & 0\\
0 & -r_1 \end{array} \right) )
+ \tau( \left( \begin{array}{cc} r_2 & 0\\
0 & -r_2 \end{array} \right) )
\\
&&\\
& = & \left( \begin{array}{cc} r_1' & 0\\ 0 & -r_1' \end{array} \right)
+ \left( \begin{array}{cc} r_2' & 0\\ 0 & -r_2' \end{array} \right)
\\
&&\\
& = &
\left( \begin{array}{cc} r_1'+r_2' & 0\\ 0 & -(r_1'+r_2') \end{array} \right)
\end{array}$$
for $r_1,r_2 \in R$,
we obtain $\eta_\tau(r_1+r_2) = (r_1 + r_2)' =
r_1'+r_2' = \eta_\tau(r_1)+\eta_\tau(r_2)$.
If $s \in F$ and $r \in R$, then
$$\begin{array}{lll}
\tau(
\left( \begin{array}{cc} sr & 0\\ 0 & -sr \end{array} \right) )
& = &
\tau( s 
\left( \begin{array}{cc} r & 0\\ 0 & -r \end{array} \right) )
\\
&&\\
& = &
s \tau(
\left( \begin{array}{cc} r & 0\\ 0 & -r \end{array} \right) )
\\
&&\\
& = & s \left( \begin{array}{cc} r' & 0\\ 0 & -r' \end{array} \right)
\\
&&\\
& = & \left( \begin{array}{cc} sr' & 0\\ 0 & -sr' \end{array} \right) \ .
\end{array}
$$
Therefore,
$$\eta_\tau(sr) = (sr)' = sr' = s \eta_\tau(r).$$
We also have
$$\begin{array}{lll}
\tau(
\left( \begin{array}{cc} 0 & r\\ 0 & 0 \end{array} \right) )
& = &
\tau( \left[
\left( \begin{array}{cc} \displaystyle{ \frac{r}{2} } & 0\\
&\\
0 & - \displaystyle{ \frac{r}{2} } \end{array} \right) ,
\left( \begin{array}{cc} 0 & 1\\ 0 & 0 \end{array} \right)
\right] )
\\
&&\\
& = &
\left[ \tau(
\left( \begin{array}{cc} \displaystyle{ \frac{r}{2} } & 0\\
&\\
0 & - \displaystyle{ \frac{r}{2} } \end{array} \right) ) ,
\tau(
\left( \begin{array}{cc} 0 & 1\\ 0 & 0 \end{array} \right) )
\right]
\\
&&\\
& = &
\left[
\left( \begin{array}{cc} \displaystyle{ \frac{r'}{2} } & 0\\
&\\
0 & - \displaystyle{ \frac{r'}{2} } \end{array} \right) ,
\left( \begin{array}{cc} 0 & 1\\ 0 & 0 \end{array} \right)
\right]
\\
&&\\
& = &
\left( \begin{array}{cc} 0 & r'\\ 0 & 0 \end{array} \right)
\end{array}$$
for $r \in r$.
Similarly we see
$$\tau(
\left( \begin{array}{cc} 0 & 0\\ r & 0 \end{array} \right) )
= \left( \begin{array}{cc} 0 & 0\\ r' & 0 \end{array} \right)
$$
for $r \in R$.
Furthermore, for $r_1,r_2 \in R$, we find
$$\begin{array}{lll}
\left( \begin{array}{cc} 0 & (r_1r_2)'\\ 0 & 0 \end{array} \right)
& = &
\tau(
\left( \begin{array}{cc} 0 & r_1r_2\\ 0 & 0 \end{array} \right) )
\\
&&\\
& = &
\tau(
\left[
\left( \begin{array}{cc} \displaystyle{ \frac{r_1}{2} } & 0\\
\\
0 & - \displaystyle{ \frac{r_1}{2} } \end{array} \right) ,
\left( \begin{array}{cc} 0 & r_2\\ 0 & 0 \end{array} \right)
\right]
)
\\
&&\\
& = &
\left[
\tau(
\left( \begin{array}{cc} \displaystyle{ \frac{r_1}{2} } & 0\\
\\
0 & - \displaystyle{ \frac{r_1}{2} } \end{array} \right) ) ,
\tau(
\left( \begin{array}{cc} 0 & r_2\\ 0 & 0 \end{array} \right) )
\right]
\\
&&\\
& = &
\left[
\left( \begin{array}{cc} \displaystyle{ \frac{r_1'}{2} } & 0\\
\\
0 & - \displaystyle{ \frac{r_1'}{2} } \end{array} \right) ,
\left( \begin{array}{cc} 0 & r_2'\\ 0 & 0 \end{array} \right)
\right]
\\
&&\\
& = &
\left( \begin{array}{cc} 0 & r_1'r_2'\\ 0 & 0 \end{array} \right) \ ,
\end{array}$$
which implies
$$\eta_\tau(r_1r_2) = (r_1r_2)' = r_1'r_2' = 
\eta_\tau(r_1) \eta_\tau(r_2).$$
Hence, $\eta_\tau \in Aut_F(R)$. Therefore,
$$Aut_F(\mathfrak{sl}_2(R)) = PGL_2(R) \cdot  Aut_F(R).$$
Let $\sigma_0 \in PGL_2(R) \cap Aut_F(R)$. Suppose that
$\sigma_0$ corresponds to conjugation by
$$\left(
\begin{array}{cc}
a & b\\ c & d \end{array}
\right)
$$
for some $a,b,c,d \in R$ with $ad - bc \in R^\times$.
Then,
$$\sigma_0(
\left(
\begin{array}{cc}
0 & 1\\ 0 & 0 \end{array}
\right) )
=
\frac{1}{ad-bc}
\left(
\begin{array}{cc}
-ac & a^2\\ -c^2 & ac \end{array}
\right)
=
\left(
\begin{array}{cc}
0 & 1\\ 0 & 0 \end{array}
\right) \ ,
$$
which implies $a = d$ and $c = 0$.
Similarly we have
$$\sigma_0(
\left(
\begin{array}{cc}
0 & 0\\ 1 & 0 \end{array}
\right) )
=
\frac{1}{ad-bc}
\left(
\begin{array}{cc}
bd & -b^2\\ d^2 & -bd \end{array}
\right)
=
\left(
\begin{array}{cc}
0 & 0\\ 1 & 0 \end{array}
\right) \ ,
$$
which implies $a = d$ and $b = 0$.
Therefore, $a = d \in R^\times$ and $b = c = 0$.
Hence, $\sigma_0$  corresponds to conjugation by
$$\left(
\begin{array}{cc}
a & 0\\ 0 & a \end{array}
\right)
$$
with $a \in R^\times$, which means $\sigma_0 = 1$.
Thus, $PGL_2(R) \cap Aut_F(R) = 1$. Since $Aut_F(R)$ normalizes
$PGL_2(R)$, we obtain
$Aut_F(\mathfrak{sl}_2(R)) = PGL_2(R) \rtimes Aut_F(R)$.
Q.E.D.
\bigskip

   Our next result is a conjugacy theorem for $s\ell_2$-triples 
in $\mathfrak{sl}_2 (R)$. We show they are all conjugate to the standard one.
In fact, we can do this for a general integral domain (which is still an $F$-algebra)
but can pick the conjugating automorphism in $PGL_2(R)$ when $R$ is a UFD.

\bigskip
\noindent
{\bf Theorem 6}.  Let the $F$ algebra $R$ be an integral domain and let  $(X,H,Y)$ be an $s\ell_2$-triple with
$X,H,Y \in \mathfrak{sl}_2(R)$.
Then there is an automorphism $\tau \in Aut_F(\mathfrak{sl}_2(R))$
satisfying
$$\begin{array}{lll}
\tau(X)
& = & \left( \begin{array}{cc} 0 & 1\\ 0 & 0 \end{array} \right) \ ,\\
&&\\
\tau(H)
& = & \left( \begin{array}{cc} 1 & 0\\ 0 & -1 \end{array} \right) \ ,\\
&&\\
\tau(Y)
& = & \left( \begin{array}{cc} 0 & 0\\ 1 & 0 \end{array} \right) \ .
\end{array}$$
Furthermore, if $R$ is a unique factorization domain then
we can choose the above $\tau$ to be in
$PGL_2(R) \subset Aut_F(\mathfrak{sl}_2(R))$
\bigskip

\noindent
{\it Proof}\quad
Suppose first  that $R$ is just an integral domain.
Put ${\mathfrak a} = FX \oplus FH \oplus FY$, and write
$$X = \left(
\begin{array}{cc}
X_1 & X_2\\ X_3 & -X_1 \end{array}
\right) \ ,
\qquad
H = \left(
\begin{array}{cc}
H_1 & H_2\\ H_3 & -H_1 \end{array}
\right) \ ,
\qquad
Y =
\left(
\begin{array}{cc}
Y_1 & Y_2\\ Y_3 & -Y_1 \end{array}
\right) \ .
$$
We define a map
$$\theta : \mathfrak{sl}_2(R) \longrightarrow \mathfrak{sl}_2(R)$$
by
$$\left(
\begin{array}{c}
Z_1'\\ Z_2'\\ Z_3' \end{array}
\right)
= B \left(
\begin{array}{c}
Z_1\\ Z_2\\ Z_3 \end{array}
\right)$$
with
$$B = \left(
\begin{array}{ccc}
H_1 & X_1 & Y_1\\ H_2 & X_2 & Y_2\\ H_3 & X_3 & Y_3
\end{array}
\right) \ ,$$
where we are using for notation the following,
$$\left(
\begin{array}{c}
Z_1\\ Z_2\\ Z_3 \end{array}
\right) = \left(
\begin{array}{cc}
Z_1 & Z_2\\ Z_3 & -Z_1\end{array}
\right) \in \mathfrak{sl}_2(R)$$
and
$$\left(
\begin{array}{c}
Z_1'\\ Z_2'\\ Z_3' \end{array}
\right) = Z_1'H + Z_2'X + Z_3'Y.$$
This map $\theta$ can clearly be extended to be
a $K$-linear map, again called $\theta$, of
$\mathfrak{sl}_2(K)$ into $\mathfrak{sl}_2(K)$. Considering
the eigenvalues of ${\rm ad}\ H$,
we see that
$\{ X,H,Y \}$ is linearly independent over $K$, which
shows $${\rm Im}\ \theta = KX \oplus KH \oplus KY =
\mathfrak{sl}_2(K).$$
Therefore, $\theta$ is an isomorphism of
$\mathfrak{sl}_2(K)$ onto $\mathfrak{sl}_2(K)$ as Lie algebras over $K$.
Since ${\rm det}\ B = H_1^2 + H_2H_3 = 1$,( use (C) in the proof of Lemma 1 for this)
we see that $B$ belongs to $GL_3(R)$, which implies
$$\mathfrak{sl}_2(R) = \theta^{-1}(RX \oplus RH \oplus RY) 
\subset \theta^{-1}(\mathfrak{sl}_2(R)) \subset \mathfrak{sl}_2(R).$$
This means $RX \oplus RH \oplus RY = \mathfrak{sl}_2(R)$.
Hence, 
$\theta$ gives an automorphism, also denoted by $\theta$,
of $\mathfrak{sl}_2(R)$ viewed as a Lie algebra over $F$.
Then, $\tau = \theta^{-1}$ is our desired conjugating automorphism.
\bigskip

\noindent
Now we suppose that $R$ is a UFD.
Let $(X,H,Y)$ be an $s\ell_2$-triple in
$\mathfrak{sl}_2(R)$.
We will follow the line of proof in Theorem 5.
By Lemma 1, we have
$$H^2 = \left( \begin{array}{cc} 1 & 0\\ 0 & 1 \end{array} \right) \ .
$$
Then, by Theorem 3,
there is $\tau' \in PGL_2(R)$ such that
$$\tau'(H) =
\left( \begin{array}{cc} 1 & 0\\ 0 & -1 \end{array} \right) \ .$$
Just as in the proof of Theorem 5,
we obtain
$$\begin{array}{lll}
\tau'(X) & = &
\left( \begin{array}{cc} 0 & u\\ 0 & 0 \end{array} \right) \\
&&\\
\tau'(Y) & = &
\left( \begin{array}{cc} 0 & 0\\ u^{-1} & 0 \end{array} \right)
\end{array}
$$
for some $u \in R^\times$. We choose $\tau'' \in PGL_2(R)$
corresponding to conjugation by
$$\left( \begin{array}{cc} u^{-1} & 0\\ 0 & 1 \end{array} \right) \ .
$$
Put $\tau = \tau'' \tau'$. Then, we have
$$\begin{array}{lll}
\tau(X)
& = & \left( \begin{array}{cc} 0 & 1\\ 0 & 0 \end{array} \right) \ ,\\
&&\\
\tau(H)
& = & \left( \begin{array}{cc} 1 & 0\\ 0 & -1 \end{array} \right) \ ,\\
&&\\
\tau(Y)
& = & \left( \begin{array}{cc} 0 & 0\\ 1 & 0 \end{array} \right) \ ,
\end{array}$$
which is as desired. Q.E.D.
\bigskip

\noindent
Proposition 2 says that we can find an $s\ell_2$-triple 
from an $H$ satisfying $H^2=I$. In the next result, we show
can also find an $s\ell_2$-triple from an  $X$ having a certain property.
We think of this as a type of Jacobson-Morozov theorem
for $\mathfrak{sl}_2$ over an integral domain. For notation we let $(r)$ denote the principal ideal 
of $R$ generated by $r \in R$. We will need the following Remark, whose proof is an easy exercise, in the proof of this result.
\bigskip

\noindent
{\bf Remark}
Suppose that $R$ is just a commutative associative $F$-algebra with identity. Then 
there is a canonical one to one correspondence between
the ideal lattice, 
${\mathcal L}(\mathfrak{sl}_2(R))$,
of the Lie algebra $\mathfrak{sl}_2(R)$ over $F$ and
the ideal lattice,
${\mathcal L}(R)$, of $R$, which is explicitly given by
the following correspondence.
$$\begin{array}{ccc}
{\mathcal L}(\mathfrak{sl}_2(R)) & \longleftrightarrow &
{\mathcal L}(R)\\
&&\\
{\mathfrak a} & \longrightarrow &
\left\{ \ a \in R\ \left| \ 
\left( \begin{array}{cc} 0&a\\ 0&0 \end{array} \right) \in {\mathfrak a}
\ \right. \right\} \\
&&\\
{\mathfrak A} \left( \begin{array}{cc} 0&1\\ 0&0 \end{array} \right)
\oplus
{\mathfrak A} \left( \begin{array}{cc} 1&0\\ 0&-1 \end{array} \right)
\oplus
{\mathfrak A} \left( \begin{array}{cc} 0&0\\ 1&0 \end{array} \right)
& \longleftarrow & {\mathfrak A}
\end{array}$$

\noindent
{\bf Theorem 7}. Assume that 
$R$ is an integral domain. Let $X \in \mathfrak{sl}_2(R)$ be a nonzero element.
Then, the following three conditions are equivalent.\par
\noindent
(1) ${\rm ad}\ X\ \mbox{is nilpotent}$ (i.e.\ $X_1^2 + X_2X_3 = 0$)
and $(X_2) + (X_3) = R$.\par
\noindent
(2) There exist $H, Y \in \mathfrak{sl}_2(R)$ such that
$(X,H,Y)$ is an $s\ell_2$-triple.\par
\noindent
(3) There is an automorphism $\rho \in Aut_F(\mathfrak{sl}_2(R))$
satisfying
$$\rho(X) =
\left(
\begin{array}{cc} 0 & 1\\ 0 & 0 \end{array}
\right) \ .$$
\bigskip

\noindent
{\it Proof}\quad

\noindent
(1) $\Rightarrow$ (2):
For each
$$X = \left(
\begin{array}{cc}
X_1 & X_2\\ X_3 & -X_1
\end{array}
\right) \ ,
$$
we see
$$({\rm ad}\ X)^3 = 4(X_1^2 + X_2X_3)({\rm ad}\ X).$$
Hence, since $X$ is not zero, we obtain that ${\rm ad}\ X$ is nilpotent
if and only if $X_1^2 + X_2X_3 = 0$, which is also
equivalent to $X^2 = O$.
Since $(X_2) + (X_3) = R$, there are $r,s \in R$ such that
$rX_2 + sX_3 = 1$.
If we put
$$H = \left(
\begin{array}{cc}
rX_2-sX_3 & 2sX_1\\ -2rX_1 & -(rX_2-sX_3)
\end{array} \right)
\quad
\mbox{and}\quad
Y = \left(
\begin{array}{cc}
-rsX_1 & s^2X_3\\ r^2X_2 & rsX_1
\end{array} \right) \ ,
$$
then we easily see
$$[H,X] = 2X,\quad [H,Y] = -2Y,\quad [X,Y] = H.$$
Therefore, we obtain that
$(X,H,Y)$ is an $s\ell_2$-triple. 
\bigskip

\noindent
(2) $\Rightarrow$ (3): This part follows from Theorem 6.
\bigskip

\noindent
(3) $\Rightarrow$ (1): Since
$${\rm ad}\ \left(
\begin{array}{cc}
0 & 1\\ 0 & 0
\end{array} \right)
$$
is nilpotent, ${\rm ad}\ X$ with
$$X = \left(
\begin{array}{cc}
X_1 & X_2\\ X_3 & -X_1
\end{array} \right)
$$
is nilpotent. Therefore, we obtain $X_1^2 + X_2X_3 = 0$.
Assume $(X_2) + (X_3) \not= R$. Then, there is
a maximal ideal ${\mathfrak M}$ such that
$$(X_2 ) + (X_3) \subset {\mathfrak M} \not= R.$$
Since $X_1^2 = -X_2X_3 \in {\mathfrak M}$, we have
$X_1 \in {\mathfrak M}$. We set
$$H = \rho^{-1}( \left(
\begin{array}{cc}
1 & 0\\ 0 & -1
\end{array} \right)
)\quad \mbox{and}\quad 
Y = \rho^{-1}( \left(
\begin{array}{cc}
0 & 0\\1 & 0
\end{array} \right)
).
$$
Because of the relations (A), (B), (C) from Lemma 1, we see that
all elements of ${\mathfrak a} = FX \oplus FH \oplus FY$
have their coefficients in ${\mathfrak M}$.
Hence, ${\mathfrak a}$ is contained in some maximal
ideal of $\mathfrak{sl}_2(R)$ by the Remark befor this Proposition. 
Since $\rho$ preserves the ideal structure of
$\mathfrak{sl}_2(R)$, we obtain that
$\mathfrak{sl}_2(F) = \rho({\mathfrak a})$
is also contained in another maximal ideal of
$\mathfrak{sl}_2(R)$, which is a contradiction.
Therefore, $(X_2) + (X_3) = R$.
This completes our proof.
Q.E.D.
\bigskip

In the case when the coordinate ring $R$ is a UFD we can 
strenghten the preceeding result by saying the conjugating 
automorphism can be taken to be in $PGL_2(R)$.
\bigskip

\noindent
{\bf Corollary 8}.  Suppose that
$R$ is a unique factorization domain. Let $X \in \mathfrak{sl}_2(R)$.
Then the following two conditions are equivalent.\par
\noindent
(1) $X_1^2 + X_2X_3 = 0$ and $(X_2) + (X_3) = R$.\par
\noindent
(2) There is $P \in GL_2(R)$ such that
$$P^{-1} X P = \left(
\begin{array}{cc}
0 & 1\\ 0 & 0
\end{array}
\right) \ .
$$
\bigskip

\noindent
{\it Proof}\quad

\noindent
(1) $\Rightarrow$ (2): Since $R$ is a unique factorization domain
and $(X_2) + (X_3) = R$,
we can write
$$\begin{array}{lll}
X_1 & = & u p_1 p_2 \cdots p_m q_1 q_2 \cdots q_n,\\
X_2 & = & v p_1^2 p_2^2 \cdots p_m^2,\\
X_3 & = & w q_1^2 q_2^2 \cdots q_n^2
\end{array}
$$
for some $u,v,w \in R^\times$ with $u^2 + vw = 0$
and for some prime
elements
$$p_1,p_2,\ldots,p_m,q_1,q_2,\ldots,q_n$$
of $R$.
Put $p_1' = up_1$.
Then, we have
$$\begin{array}{lll}
X_1 & = & (p_1' p_2 \cdots p_m)(q_1 q_2 \cdots q_n),\\
&&\\
X_2 & = &  \displaystyle{ -\frac{1}{w} }(p_1'^2 p_2^2 \cdots p_m^2),\\
&&\\
X_3 & = & w (q_1^2 q_2^2 \cdots q_n^2).
\end{array}
$$
We set $p = p_1' p_2 \cdots p_m$ and $q = q_1 q_2 \cdots q_n$.
Then,
$$X_1 = pq,\quad X_2 = - \frac{1}{w} p^2,\quad X_3 = wq^2.$$
Since $(X_2)+(X_3) = R$, we can find  elements $r,s \in R$
such that
$$rX_2 + sX_3 = 1.$$
We define $P$ by
$$P = \left(
\begin{array}{cc}
p & swq\\ wq & rp
\end{array} \right) \ .
$$
Then we have
$$\begin{array}{lll}
\left(
\begin{array}{cc}
X_1 & X_2\\ X_3 & -X_1
\end{array} \right)
P
& = &
\left(
\begin{array}{cc}
pq & \displaystyle{ -\frac{1}{w} }p^2\\
wq^2 & -pq
\end{array} \right)
\ \left(
\begin{array}{cc}
p & swq\\ wq & rp
\end{array} \right)
\\
&&\\
& = &
\left(
\begin{array}{cc}
p^2q-p^2q & swpq^2-r\displaystyle{\frac{1}{w}}p^3\\
wpq^2-wpq^2 & w^2sq^3-rp^2q
\end{array} \right)
\\
&&\\
& = &
\left(
\begin{array}{cc}
0 & p\\ 0 & wq
\end{array} \right)
\end{array}
$$
and
$$P
\left(
\begin{array}{cc}
0 & 1\\ 0 & 0
\end{array} \right)
=
\left(
\begin{array}{cc}
0 & p\\ 0 & wq
\end{array} \right) \ ,
$$
which implies
$$
X P =
P
\left(
\begin{array}{cc}
0 & 1\\ 0 & 0
\end{array} \right) \ .
$$
Since
$$\begin{array}{lll}
{\rm det}\ P & = &
rp^2 - sw^2q^2\\
&&\\
& = & (-w)(- \displaystyle{ \frac{1}{w} } rp^2 + wsq^2)\\
&&\\
& = & (-w)(rX_2 + sX_3)\\
&&\\
& = & -w,
\end{array}$$
we have
$P \in GL_2(R)$. Therefore, we obtain
$$P^{-1}XP = 
\left(
\begin{array}{cc}
0 & 1\\ 0 & 0
\end{array} \right) \ .
$$
(2) $\Rightarrow$ (1): This part follows directly from Theorem 7.
 Q.E.D.
\bigskip

\noindent
{\bf  Remark}. In the above Corollary we used Theorem 7 to deduce 
that the second statement implies the first. However we can do this directly,
even when $R$ is just a commutative $F$ algebra.We present the argument here.

Since
$$X = P
\left(
\begin{array}{cc}
0 & 1\\ 0 & 0
\end{array} \right) P^{-1}\ ,
$$
we see $X^2 = O$ and so $X_1^2 + X_2X_3 = 0$. We write
$$P =
\left(
\begin{array}{cc}
a & b\\ c & d
\end{array} \right)
$$
for some $a,b,c,d \in R$ with $ad - bc \in R^\times$.
Then we have
$$X = P
\left(
\begin{array}{cc}
0 & 1\\ 0 & 0
\end{array} \right)
P^{-1} =
\frac{1}{ad-bc}
\left(
\begin{array}{cc}
-ac & a^2\\ -c^2 & ac
\end{array} \right) \ .
$$
We set
$$t = \frac{d}{ad-bc}\quad \mbox{and}\quad u = -\frac{b}{ad-bc}\ .$$
Then $1 = ta + uc$ and $ac = ta^2c + uac^2$.
On the other hand we obtain
$$\begin{array}{lll}
1 & = & (ta+uc)^2\\
& = & t^2a^2 + 2tuac + u^2c^2\\
& = & t^2a^2 + 2tu(ta^2c + uac^2) + u^2c^2\\
& = & t^2(1+2uc)a^2 + u^2(1+2ta)c^2.
\end{array}$$
This means
$$(X_2) + (X_3) = (a^2) + (c^2) = R,$$
which is what we want.
\bigskip

\noindent
We close this section by investigateing when an element
$$Z = \left(
\begin{array}{cc}
Z_1 & Z_2\\ Z_3 & -Z_1
\end{array}
\right) \in \mathfrak{sl}_2(R)$$
belongs to a subalgebra of $\mathfrak{sl}_2(R)$ which is isomorphic to $\mathfrak{sl}_2(F)$.
 The following result will provide an answer to this question.
\bigskip

\noindent
{\bf Proposition 9}.  Suppose that
$R$ is a unique factorization domain and
that $F$ is square root closed, that is, $F = F^2$.
Let $Z$ be a nonzero element of $\mathfrak{sl}_2(R)$.
Then the following two conditions are equivalent.\par
\noindent
(1) There is a subalgebra ${\mathfrak a}$ of $\mathfrak{sl}_2(R)$
such that $Z \in {\mathfrak a}$ and
${\mathfrak a} \simeq \mathfrak{sl}_2(F)$.\par
\noindent
(2) If $Z_1^2 + Z_2Z_3 = 0$, then $(Z_2) + (Z_3) = R$, otherwise
$Z_1^2 + Z_2Z_3 \in F^\times$.
\bigskip

\noindent
{\it Proof}\quad
\noindent
(1) $\Rightarrow$ (2): By Theorem 6, we can find a suitable
$\tau \in PGL_2(R) \subset Aut_F(\mathfrak{sl}_2(R))$ such that
$\tau({\mathfrak a}) = \mathfrak{sl}_2(F)$.
Say $\tau$ is conjugation by the matrix
$$P = \left(
\begin{array}{cc}
a & b\\ c & d
\end{array} \right)
\in GL_2(R)$$
and put $u = ad - bc \in R^\times$.
By the assumption, we have
$$PZP^{-1} =
\left(
\begin{array}{cc}
Z_1' & Z_2'\\ Z_3' & -Z_1'
\end{array} \right)
\in \mathfrak{sl}_2(F)$$
with $Z_1',Z_2',Z_3' \in F$.
Then,
$$Z_1^2+Z_2Z_3 = -\ {\rm det}\ Z 
= -\ {\rm det}\ PZP^{-1} = Z_1'^2+Z_2'Z_3' \in F.$$
If $Z_1'^2 + Z_2'Z_3' \not= 0$, then we are done.
We suppose $Z_1'^2 + Z_2'Z_3' = 0$, which
implies $Z_1^2 + Z_2Z_3 = 0$.
Since
$$\left(
\begin{array}{cc}
Z_1' & Z_2'\\ Z_3' & -Z_1'
\end{array} \right)
=
\left(
\begin{array}{cc}
a & b\\ c & d
\end{array} \right)
\ \left(
\begin{array}{cc}
Z_1 & Z_2\\ Z_3 & -Z_1
\end{array} \right)
\ \left(
\begin{array}{cc}
a & b\\ c & d
\end{array} \right)^{-1}\ ,$$
we obtain the following.
$$\left\{
\begin{array}{llrcrcr}
uZ_1' & = &
(ad+bc)Z_1 & -& acZ_2 & + & bdZ_3\\
uZ_2' & = & - 2abZ_1 & + & a^2Z_2 & - & b^2Z_3\\
uZ_3' & = & 2cdZ_1 & - & c^2Z_2 & + & d^2Z_3
\end{array}
\right.
$$
We put
$$A = \left(
\begin{array}{ccc}
ad+bc & -ac & bd\\
-2ab & a^2 & -b^2\\
2cd & -c^2 & d^2
\end{array} \right) \ .
$$
Then, we can rewrite the above equations as
$$\left(
\begin{array}{c}
Z_1'\\ Z_2'\\ Z_3'
\end{array}
\right)
=
\frac{1}{u}\ 
A\ 
\left(
\begin{array}{c}
Z_1\\ Z_2\\ Z_3
\end{array}
\right) \ .
$$
Since ${\rm det}\ A = u^3$, we see $A \in GL_3(R)$.
Now we suppose $(Z_2) + (Z_3) \not= R$. Then, there is a
maximal ideal ${\mathfrak M}$ of $R$ such that
$$(Z_2) + (Z_3) \subset {\mathfrak M} \not= R.$$ 
Because $Z_1^2 = -Z_2Z_3 \in {\mathfrak M}$, we observe
$Z_1 \in {\mathfrak M}$. Therefore
$$\left(
\begin{array}{c}
Z_1\\ Z_2\\ Z_3
\end{array}
\right)
\in {\mathfrak M}^3,$$
which implies $Z_1',Z_2',Z_3' \in {\mathfrak M}$.
Since $F \cap {\mathfrak M} = (0)$ we find that
$Z_1' = Z_2' = Z_3' = 0$.
This means
$$A
\left(
\begin{array}{c}
Z_1\\ Z_2\\ Z_3
\end{array}
\right)
=
\left(
\begin{array}{c}
0\\ 0\\ 0
\end{array}
\right) \ ,
$$
which shows
$$\left(
\begin{array}{c}
Z_1\\ Z_2\\ Z_3
\end{array}
\right)
= A^{-1}
\left(
\begin{array}{c}
0\\ 0\\ 0
\end{array}
\right)
=
\left(
\begin{array}{c}
0\\ 0\\ 0
\end{array}
\right) \ .
$$
This is a contradiction.
Hence, $(Z_2) + (Z_3) = R$.
(Note that to prove this part of the result we did  not use the hypothesis that $F = F^2$.)
\bigskip

\noindent
(2) $\Rightarrow$ (1): First, we suppose that
$$Z_1^2 + Z_2Z_3 = 0\quad \mbox{and}\quad
(Z_2) + (Z_3) = R.$$
Then by Theorem 7 we have an $s\ell_2$-triple
$(Z,H,Y)$ for some $H,Y \in \mathfrak{sl}_2(R)$.
\bigskip

\noindent
Next, we
suppose that
$$Z_1^2 + Z_2Z_3 \not= 0\quad \mbox{and}\quad
Z_1^2 + Z_2Z_3 \in F^\times.$$
Since $F^\times = {F^\times}^2$, we can choose
$\lambda \in F^\times$ such that $Z_1^2 + Z_2Z_3 = \lambda^2$.
We put
$$Z' = \frac{1}{\lambda}\ Z \in \mathfrak{sl}_2(R).$$
Then, we see
$$Z'^2 = \frac{1}{\lambda^2}\ 
\left(
\begin{array}{cc}
Z_1^2+Z_2Z_3 & 0\\ 0 & Z_1^2+Z_2Z_3
\end{array}
\right)
=
\left(
\begin{array}{cc}
1 & 0\\ 0 & 1
\end{array}
\right) \ .
$$
Hence by Proposition 2 we have an $s\ell_2$-triple 
$(X,Z',Y)$ for some $X,Y \in \mathfrak{sl}_2(R)$.
Then,
$$Z \in {\mathfrak a} = FX \oplus FZ' \oplus FY \simeq
\mathfrak{sl}_2(F).$$
\bigskip

\noindent
Thus 
we have found a subalgebra ${\mathfrak a}$ of
$\mathfrak{sl}_2(R)$ such that
$Z \in {\mathfrak a}$ and ${\mathfrak a} \simeq \mathfrak{sl}_2(F)$.
Q.E.D.
\bigskip

\section{Isomorphisms and Derivations of $\mathfrak{sl}_2(R)$.} In this, our final section, we investigate two questions. The first is when can we have two Lie algebras $\mathfrak{sl}_2(R)$ and $\mathfrak{sl}_2(S)$ be isomorphic when both $R$ and $S$ are commutative algbras over $F$ and $R$ is an integral domain. We will see this happens only when $R$ and $S$ are isomorphic. Finally we go on to study the derivation algebra of $\mathfrak{sl}_2(R)$ and then make some closing remarks.
\bigskip

\noindent
{\bf Proposition 10}. Suppose that
$R$ is an integral domain and $S$ is just a commutative
$F$-algebra. Then the following two conditions are equivalent.\par
\noindent
(1) $\mathfrak{sl}_2(R) \simeq \mathfrak{sl}_2(S)$
as Lie algebras over $F$.\par
\noindent
(2) $R \simeq S$ as $F$-algebras.
\bigskip

\noindent
{\it Proof}\quad
\noindent
(1) $\Rightarrow$ (2): Let
$$\phi : \mathfrak{sl}_2(S) \longrightarrow \mathfrak{sl}_2(R)$$
be an isomorphism. We take
$${\mathfrak s} = F
\left(
\begin{array}{cc}
0 & 1\\ 0 & 0
\end{array}
\right)
\oplus F
\left(
\begin{array}{cc}
1 & 0\\ 0 & -1
\end{array}
\right)
\oplus F
\left(
\begin{array}{cc}
0 & 0\\ 1 & 0
\end{array}
\right)
= \mathfrak{sl}_2(F)
\subset \mathfrak{sl}_2(S).$$
We put ${\mathfrak a} = \phi({\mathfrak s}) \subset \mathfrak{sl}_2(R)$
and
$$\begin{array}{lll}
X & = & \phi(
\left(
\begin{array}{cc}
0 & 1\\ 0 & 0
\end{array}
\right)
),\\
&&\\
H & = & \phi(
\left(
\begin{array}{cc}
1 & 0\\ 0 & -1
\end{array}
\right)
),\\
&&\\
Y & = & \phi(
\left(
\begin{array}{cc}
0 & 0\\ 1 & 0
\end{array}
\right)
).
\end{array}$$
Then we find an automorphism $\tau \in 
Aut(\mathfrak{sl}_2(R))$ such that
$$\tau({\mathfrak a}) = \mathfrak{sl}_2(F) \subset \mathfrak{sl}_2(R)$$
by Theorem 6. Therefore, we may assume
$\phi(\mathfrak{sl}_2(F)) = \mathfrak{sl}_2(F)$
from the beginning. Again Theorem 6 implies that we may also assume
$$\begin{array}{lll}
X
& = & \left( \begin{array}{cc} 0 & 1\\ 0 & 0 \end{array} \right) \ ,\\
&&\\
H
& = & \left( \begin{array}{cc} 1 & 0\\ 0 & -1 \end{array} \right) \ ,\\
&&\\
Y
& = & \left( \begin{array}{cc} 0 & 0\\ 1 & 0 \end{array} \right) \ .
\end{array}$$
Using the eigenspaces of
$$\left( \begin{array}{cc} 1 & 0\\ 0 & -1 \end{array} \right)
$$
with eigenvalues $0,\pm 2$
in $\mathfrak{sl}_2(S)$ and $\mathfrak{sl}_2(R)$ respectively,
we have
$$\begin{array}{lll}
\phi(\ S
\left( \begin{array}{cc} 0 & 1\\ 0 & 0 \end{array} \right)
\ )
& = & R \left( \begin{array}{cc} 0 & 1\\ 0 & 0 \end{array} \right) \ ,
\\
&&\\
\phi(\ S
\left( \begin{array}{cc} 1 & 0\\ 0 & -1 \end{array} \right)
\ )
& = & R \left( \begin{array}{cc} 1 & 0\\ 0 & -1 \end{array} \right) \ ,
\\
&&\\
\phi(\ S
\left( \begin{array}{cc} 0 & 0\\ 1& 0 \end{array} \right)
\ )
& = & R \left( \begin{array}{cc} 0 & 0\\ 1 & 0 \end{array} \right) \ .
\end{array}$$
Now we define $\phi^*$ by $\phi^*(s) = r$ if
$$
\phi(
\left( \begin{array}{cc} s & 0\\ 0 & -s \end{array} \right)
)
= \left( \begin{array}{cc} r & 0\\ 0 & -r \end{array} \right)
$$
with $s \in S$ and $r \in R$. This $\phi^*$ is a bijective
map of $S$ to $R$ since $\phi$ is an isomorphism.
Then we have
$$\begin{array}{lll}
\left( \begin{array}{cc} 0 & r\\ 0 & 0 \end{array} \right)
& = &
\left[
\left( \begin{array}{cc} r & 0\\ 0 & -r \end{array} \right) ,
\left( \begin{array}{cc} 0 & \displaystyle{ \frac{1}{2} }
\\ 0 & 0 \end{array} \right)
\right]
\\
&&\\
& = & 
\left[
\phi(
\left( \begin{array}{cc} s & 0\\ 0 & -s \end{array} \right)
) ,
\phi(
\left( \begin{array}{cc} 0 & \displaystyle{ \frac{1}{2} }
\\ 0 & 0 \end{array} \right)
)
\right]
\\
&&\\
& = &
\phi(
\left[
\left( \begin{array}{cc} s & 0\\ 0 & -s \end{array} \right) ,
\left( \begin{array}{cc} 0 & \displaystyle{ \frac{1}{2} }
\\ 0 & 0 \end{array} \right)
\right]
)\\
&&\\
& = &
\phi(
\left( \begin{array}{cc} 0 & s\\ 0 & 0 \end{array} \right)
)
\end{array}$$
and
$$\phi(
\left( \begin{array}{cc} 0 & s\\ 0 & 0 \end{array} \right)
) =
\left( \begin{array}{cc} 0 & r\\ 0 & 0 \end{array} \right) \ .$$
Similarly we obtain
$$\phi(
\left( \begin{array}{cc} 0 & 0\\ s & 0 \end{array} \right)
) =
\left( \begin{array}{cc} 0 & 0\\ r & 0 \end{array} \right) \ .$$
Since $\phi$ is additive,
we see $\phi^*(s_1+s_2) = \phi^*(s_1) + \phi^*(s_2)$
for $s_1,s_2 \in S$. We also have
$$
\begin{array}{lll}
\phi(
\left( \begin{array}{cc} s_1s_2 & 0\\ 0 & -s_1s_2 \end{array} \right)
) & = &
\phi(
[
\left( \begin{array}{cc} 0 & s_1\\ 0 & 0 \end{array} \right) ,
\left( \begin{array}{cc} 0 & 0\\ s_2 & 0 \end{array} \right)
]
)
\\
&&\\
& = &
[
\phi(
\left( \begin{array}{cc} 0 & s_1\\ 0 & 0 \end{array} \right)
) , \phi(
\left( \begin{array}{cc} 0 & 0\\ s_2 & 0 \end{array} \right)
)
]
\\
&&\\
& = &
[
\left( \begin{array}{cc} 0 & r_1\\ 0 & 0 \end{array} \right) ,
\left( \begin{array}{cc} 0 & 0\\ r_2 & 0 \end{array} \right)
]
\\
&&\\
& = & 
\left( \begin{array}{cc} r_1r_2 & 0
\\ 0 & -r_1r_2 \end{array} \right)
\end{array}$$
if $\phi^*(s_i) = r_i$ for $i = 1,2$.
Hence, we obtain
$$\phi^*(s_1s_2) = \phi^*(s_1) \phi^*(s_2).$$
We already know that $\phi^*(1) = 1$ as well as that
$\phi^*$ is $F$-linear.
Therefore we have shown that $\phi^*$ is an $F$-algebra isomorphism of
$S$ onto $R$.
\bigskip

\noindent
Since  (2) $\Rightarrow$ (1) is trivial we are done.
Q.E.D.
\bigskip

Next, we will determine the derivations of $\mathfrak{sl}_2(R)$.
Here, $Der_F$ gives all $F$-derivations, and $InnD$ gives all inner
derivasions. One can view $InnD(\mathfrak{sl}_2(R))$ and
$Der_F(R)$ as subalgebras in the full derivation Lie algebra
$Der_F(\mathfrak{sl}_2(R))$ of $\mathfrak{sl}_2(R)$. Indeed, it is well known that 
$InnD(\mathfrak{sl}_2(R))$ is an ideal of $Der_F(\mathfrak{sl}_2(R))$. Also, $Der_F(R)$ acts on
$\mathfrak{sl}_2(R)$ by saying for $D \in Der_F(R)$ we have
$$  D(\left(
\begin{array}{cc}
a & b\\ c & -a
\end{array}
\right)) =\left(
\begin{array}{cc}
D(a)& D(b)\\D(c)&-D(a)
\end{array}.
\right)$$
It is clear that this defines a Lie algebra injection of $Der_F(R)$ into $Der_F(\mathfrak{sl}_2(R))$.
We thus use this to identify $Der_F(R)$ as a subalgebra of $Der_F(\mathfrak{sl}_2(R))$. Below we show that 
$Der_F(\mathfrak{sl}_2(R))$ is the semi-direct product of the ideal $InnD(\mathfrak{sl}_2(R))$ and the subalgebra $Der_F(R)$. We can do this in a general setting so only need to assume that $R$ is a commutative $F$-algebra.
\bigskip

\noindent
{\bf Proposition 11}.  Suppose that
$R$ is a commutative $F$-algebra. Then,
$$Der_F(\mathfrak{sl}_2(R)) = InnD(\mathfrak{sl}_2(R))  \oplus\ Der_F(R).$$
\bigskip

\noindent
{\it Proof}\quad
In our proof we will use the following easily established relations:
$$\begin{array}{lll}
{[}
\left( \begin{array}{cc} 1 & 0\\ 0 & -1 \end{array} \right) ,
\left( \begin{array}{cc} a & b\\ c & -a \end{array} \right)
{]}
& = &
\left( \begin{array}{cc} 0 & 2b\\ -2c & 0 \end{array} \right) \ ,\\
&&\\
{[}
\left( \begin{array}{cc} 0 & 1\\ 0 & 0 \end{array} \right) ,
\left( \begin{array}{cc} a & b\\ c & -a \end{array} \right)
{]}
& = &
\left( \begin{array}{cc} c & -2a\\ 0 & -c \end{array} \right) \ ,\\
&&\\
{[}
\left( \begin{array}{cc} 0 & 0\\ 1 & 0 \end{array} \right) ,
\left( \begin{array}{cc} a & b\\ c & -a \end{array} \right)
{]}
& = &
\left( \begin{array}{cc} -b & 0\\ 2a & b \end{array} \right) \ .
\end{array}$$
Let $D \in Der_F(\mathfrak{sl}_2(R))$.
If we write
$$D(
\left( \begin{array}{cc} 0 & 1\\ 0 & 0 \end{array} \right)
) = 
\left( \begin{array}{cc} a & b\\ c & -a \end{array} \right)
$$
and
$$D(
\left( \begin{array}{cc} 0 & 0\\ 1 & 0 \end{array} \right)
) =
\left( \begin{array}{cc} e & f\\ g & -e \end{array} \right)
$$
for some $a,b,c,e,f,g \in R$, then we have
$$\begin{array}{lll}
D(
\left( \begin{array}{cc} 1 & 0\\ 0 & -1 \end{array} \right)
) & = & 
D(
{[}
\left( \begin{array}{cc} 0 & 1\\ 0 & 0 \end{array} \right) \ ,
\left( \begin{array}{cc} 0 & 0\\ 1 & 0 \end{array} \right)
{]}
)\\
&&\\
& = &
{[}
D(
\left( \begin{array}{cc} 0 & 1\\ 0 & 0 \end{array} \right)
) ,
\left( \begin{array}{cc} 0 & 0\\ 1 & 0 \end{array} \right)
{]}
+
{[}
\left( \begin{array}{cc} 0 & 1\\ 0 & 0 \end{array} \right) \ ,
D(
\left( \begin{array}{cc} 0 & 0\\ 1 & 0 \end{array} \right)
)
{]}
\\
&&\\
& = &
\left( \begin{array}{cc} b & 0\\ -2a & -b \end{array} \right)
+
\left( \begin{array}{cc} g & -2e\\ 0 & -g \end{array} \right)
\\
&&\\
& = & \left( \begin{array}{cc} b+g & -2e\\ -2a & -b-g \end{array} \right) \ .
\end{array}$$
Since
$$\begin{array}{ll}
& 2 \left( \begin{array}{cc} a & b\\ c & -a \end{array} \right)
\\
&\\
= &
2\ D(
\left( \begin{array}{cc} 0 & 1\\ 0 & 0 \end{array} \right)
)\\
&\\
= &
D(
\left( \begin{array}{cc} 0 & 2\\ 0 & 0 \end{array} \right)
)\\
&\\
= &
D(
{[}
\left( \begin{array}{cc} 1 & 0\\ 0 & -1 \end{array} \right) ,
\left( \begin{array}{cc} 0 & 1\\ 0 & 0 \end{array} \right)
{]}
)\\
&\\
= &
{[}
D(
\left( \begin{array}{cc} 1 & 0\\ 0 & -1 \end{array} \right)
) ,
\left( \begin{array}{cc} 0 & 1\\ 0 & 0 \end{array} \right) ,
{]}
+
{[}
\left( \begin{array}{cc} 1 & 0\\ 0 & -1 \end{array} \right) ,
D(
\left( \begin{array}{cc} 0 & 1\\ 0 & 0 \end{array} \right)
)
{]}
\\
&\\
= & 
{[}
\left( \begin{array}{cc} b+g & -2e\\ -2a & -b-g \end{array} \right) ,
\left( \begin{array}{cc} 0 & 1\\ 0 & 0 \end{array} \right)
{]}
+
{[}
\left( \begin{array}{cc} 1 & 0\\ 0 & -1 \end{array} \right) ,
\left( \begin{array}{cc} a & b\\ c & -a \end{array} \right)
{]}
\\
&\\
= &
\left( \begin{array}{cc} 2a & 2(b+g)\\ 0 & -2a \end{array} \right)
+
\left( \begin{array}{cc} 0 & 2b\\ -2c & 0 \end{array} \right)
\\
&\\
= &
2 \left( \begin{array}{cc} a & 2b+g\\ -c & -a \end{array} \right) \ ,
\end{array}$$
we obtain $c = 0$ and $b + g = 0$. Similarly we have
$$\begin{array}{ll}
& - 2 \left( \begin{array}{cc} e & f\\ g & -e \end{array} \right)
\\
&\\
= &
-2\ D(
\left( \begin{array}{cc} 0 & 0\\ 1 & 0 \end{array} \right)
)\\
&\\
= &
D(
\left( \begin{array}{cc} 0 & 0\\ -2 & 0 \end{array} \right)
)\\
&\\
= &
D(
{[}
\left( \begin{array}{cc} 1 & 0\\ 0 & -1 \end{array} \right) ,
\left( \begin{array}{cc} 0 & 0\\ 1 & 0 \end{array} \right)
{]}
)\\
&\\
= &
{[}
D(
\left( \begin{array}{cc} 1 & 0\\ 0 & -1 \end{array} \right)
) ,
\left( \begin{array}{cc} 0 & 0\\ 1 & 0 \end{array} \right) ,
{]}
+
{[}
\left( \begin{array}{cc} 1 & 0\\ 0 & -1 \end{array} \right) ,
D(
\left( \begin{array}{cc} 0 & 0\\ 1 & 0 \end{array} \right)
)
{]}
\\
&\\
= &
{[}
\left( \begin{array}{cc} b+g & -2e\\ -2a & -b-g \end{array} \right) ,
\left( \begin{array}{cc} 0 & 0\\ 1 & 0 \end{array} \right)
{]}
+
{[}
\left( \begin{array}{cc} 1 & 0\\ 0 & -1 \end{array} \right) ,
\left( \begin{array}{cc} e & f\\ g & -e \end{array} \right)
{]}
\\
&\\
= &
\left( \begin{array}{cc} -2e & 0\\ -2(b+g) & 2e \end{array} \right)
+
\left( \begin{array}{cc} 0 & 2f\\ -2g & 0 \end{array} \right)
\\
&\\
= &
-2 \left( \begin{array}{cc} e & -f\\ b+2g & -e \end{array} \right) \ ,
\end{array}$$
which implies $f = 0$ and $b + g = 0$.
Therefore, we have 
$$D\ :\ \left\{ \ 
\begin{array}{lll}
\left( \begin{array}{cc} 0 & 1\\ 0 & 0 \end{array} \right)
& \mapsto &
\left( \begin{array}{cc} a & b\\ 0 & -a \end{array} \right) \ ,
\\
&&\\
\left( \begin{array}{cc} 1 & 0\\ 0 & -1 \end{array} \right)
& \mapsto &
\left( \begin{array}{cc} 0 & -2e\\ -2a & 0 \end{array} \right) \ ,
\\
&&\\
\left( \begin{array}{cc} 0 & 0\\ 1 & 0 \end{array} \right)
& \mapsto &
\left( \begin{array}{cc} e & 0\\ -b & -e \end{array} \right) \ .
\end{array}
\right.
$$
Hence, we find
$$D = {\rm ad}\ 
\left( \begin{array}{cc} \displaystyle{ \frac{b}{2} }
& e\\ -a &
\displaystyle{ - \frac{b}{2} } \end{array} \right)
$$
on $\mathfrak{sl}_2(F)$,
which means
$$\left. (D - {\rm ad}\ 
\left( \begin{array}{cc} \displaystyle{ \frac{b}{2} }
& e\\ -a &
\displaystyle{ - \frac{b}{2} } \end{array} \right)
)
\right|_{\mathfrak{sl}_2(F)} = 0.$$
We put
$$D' = D - {\rm ad}\ 
\left( \begin{array}{cc} \displaystyle{ \frac{b}{2} }
& e\\ -a &
\displaystyle{ - \frac{b}{2} } \end{array} \right) \ .$$
Now, for each $r \in R$ we write
$$D'(
\left( \begin{array}{cc} r & 0\\ 0 & -r \end{array} \right)
) =
\left( \begin{array}{cc} E_{11}(r) & E_{12}(r)\\
E_{21}(r) & -E_{11}(r) \end{array} \right)
$$
with $E_{ij}(r) \in R$ for $(i,j) = (1,1), (1,2), (2,1)$.
Since
$$\begin{array}{lll}
D'(
\left( \begin{array}{cc} 0 & s\\ 0 & 0 \end{array} \right)
) & = &
D'(
\left[
\left( \begin{array}{cc} s & 0\\ 0 & -s \end{array} \right) ,
\left( \begin{array}{cc} 0 & \displaystyle{ \frac{1}{2} }
\\ 0 & 0 \end{array} \right)
\right]
)
\\
&&\\
& = &
\left[
\left( \begin{array}{cc} E_{11}(s) & E_{12}(s)\\
E_{21}(s) & -E_{11}(s) \end{array} \right) ,
\left( \begin{array}{cc} 0 & \displaystyle{ \frac{1}{2} }\\
0 & 0 \end{array} \right)
\right]
\\
&&\\
& = &
\displaystyle{ \frac{1}{2} }\ 
\left( \begin{array}{cc} -E_{21}(s) & 2E_{11}(s)\\
0 & -E_{21}(s) \end{array} \right)
\end{array}$$
for each $s \in R$, we have
$$\begin{array}{ll}
& \left( \begin{array}{cc} -E_{21}(rs) & 2E_{11}(rs)\\
0 & E_{21}(rs) \end{array} \right)
\\
&\\
= &
\displaystyle{ \frac{1}{2} }\ 
\left( \begin{array}{cc} -E_{21}(2rs) & 2E_{11}(2rs)\\
0 & E_{21}(2rs) \end{array} \right)
\\
&\\
= &
D'(
\left( \begin{array}{cc} 0 & 2rs\\ 0 & 0 \end{array} \right)
)
\\
&\\
= &
D'(
{[}
\left( \begin{array}{cc} r & 0\\ 0 & -r \end{array} \right) ,
\left( \begin{array}{cc} 0 & s\\ 0 & 0 \end{array} \right)
{]}
)
\\
&\\
= &
{[}
\left( \begin{array}{cc} E_{11}(r) & E_{12}(r)\\
E_{21}(r) & -E_{11}(r) \end{array} \right) ,
\left( \begin{array}{cc} 0 & s\\ 0 & 0 \end{array} \right)
{]}
\\
&\\
& \qquad
+
{[}
\left( \begin{array}{cc} r & 0\\ 0 & -r \end{array} \right) ,
\displaystyle{ \frac{1}{2} }
\left( \begin{array}{cc} -E_{21}(s) & 2E_{11}(s)\\
0 & E_{21}(s) \end{array} \right)
{]}
\\
&\\
= &
\left( \begin{array}{cc} -E_{21}(r)s & 2E_{11}(r)s\\
0 & E_{21}(r)s \end{array} \right)
+
\left( \begin{array}{cc} 0 & 2rE_{11}(s)\\
0 & 0 \end{array} \right)
\\
&\\
= &
\left( \begin{array}{cc} -E_{21}(r)s & 2E_{11}(r)s + 2r E_{11}(s)\\
0 & E_{21}(r)s \end{array} \right)
\end{array}$$
for all $r,s \in R$. Thus, we obtain
$$\begin{array}{l}
E_{21}(rs) = E_{21}(r)s\\
E_{11}(rs) = E_{11}(r)s + r E_{11}(s)
\end{array}$$
for all $r,s \in R$. Immediately we see $E_{11} \in Der_F(R)$.
On the other hand, by the fact that
$$D'(
\left( \begin{array}{cc} 0 & 0\\ 1 & 0 \end{array} \right)
) =
\left( \begin{array}{cc} 0 & 0\\ 0 & 0 \end{array} \right)
$$
implies $E_{21}(1) = 0$, we have
$$\begin{array}{lll}
E_{21}(s) & = & E_{21}(1\cdot s)\\
& = & E_{21}(1)s\\
& = & 0\cdot s\\
& = & 0
\end{array}$$
for all $s \in R$.
This means that $E_{21} = 0$.
Then,
$$\begin{array}{lll}
D'(
\left( \begin{array}{cc} r & 0\\ 0 & -r \end{array} \right)
) & = &
\left( \begin{array}{cc} E_{11}(r) & E_{12}(r)\\
0 & -E_{11}(r) \end{array} \right) \ ,
\\
&&\\
D'(
\left( \begin{array}{cc} 0 & s\\ 0 & 0 \end{array} \right)
)
& = &
\left( \begin{array}{cc} 0 & E_{11}(s)\\ 0 & 0 \end{array} \right)
\end{array}
$$
for all $r,s \in R$.
Similarly we  obtain $E_{12} = 0$ as well as
$$D'(
\left( \begin{array}{cc} 0 & 0\\ s & 0 \end{array} \right)
)
=
\left( \begin{array}{cc} 0 & 0\\ E_{11}(s) & 0 \end{array} \right)
$$
for all $s \in R$.
In particular,
$$D'(
\left( \begin{array}{cc} r & 0\\ 0 & -r \end{array} \right)
) =
\left( \begin{array}{cc} E_{11}(r) & 0\\
0 & -E_{11}(r) \end{array} \right) \ .
$$
Therefore, we have
$$\begin{array}{lll}
D'(
\left( \begin{array}{cc} a & b\\ c & -a \end{array} \right)
) & = &
\left( \begin{array}{cc} E_{11}(a) & E_{11}(b)\\
E_{11}(c) & -E_{11}(a) \end{array} \right)
\\
&&\\
& = &
E_{11}(
\left( \begin{array}{cc} a & b\\ c & -a \end{array} \right)
)
\end{array}
$$
for all
$$\left( \begin{array}{cc} a & b\\ c & -a \end{array} \right)
\in \mathfrak{sl}_2(R).$$
This shows $D' = E_{11} \in Der_F(R)$.
Hence, we have
$$Der_F(\mathfrak{sl}_2(R)) = InnD(\mathfrak{sl}(R)) + Der_F(R).$$
We take
$$D'' \in InnD(\mathfrak{sl}_2(R)) \cap Der_F(R).$$
Then, $D'' = {\rm ad}\ Z$ for some
$$Z =
\left( \begin{array}{cc} Z_1 & Z_2\\ Z_3 & -Z_1 \end{array} \right)
\in \mathfrak{sl}_2(R).
$$
We note that $D''(1) = 0$, since
$D''(1) = D(1 \cdot 1) = D''(1) + D''(1)$.
Therefore, we have
$$\begin{array}{lll}
\left( \begin{array}{cc} 0 & 0\\ 0 & 0 \end{array} \right)
& = &
\left( \begin{array}{cc} 0 & D''(1)\\ 0 & 0 \end{array} \right)
\\
&&\\
& = &
D''(
\left( \begin{array}{cc} 0 & 1\\ 0 & 0 \end{array} \right)
)
\\
&&\\
& = &
{[}
Z ,
\left( \begin{array}{cc} 0 & 1\\ 0 & 0 \end{array} \right)
{]}
\\
&&\\
& = &
\left( \begin{array}{cc} -Z_3 & 2Z_1\\ 0 & Z_3 \end{array} \right)
\end{array}$$
and $Z_1 = Z_3 = 0$.
Similarly, we obtain
$$\begin{array}{lll}
\left( \begin{array}{cc} 0 & 0\\ 0 & 0 \end{array} \right)
& = &
\left( \begin{array}{cc} 0 & 0\\ D''(1) & 0 \end{array} \right)
\\
&&\\
& = &
D''(
\left( \begin{array}{cc} 0 & 0\\ 1 & 0 \end{array} \right)
)
\\
&&\\
& = &
{[}
Z ,
\left( \begin{array}{cc} 0 & 0\\ 1 & 0 \end{array} \right)
{]}
\\
&&\\
& = &
\left( \begin{array}{cc} Z_2 & 0\\ -2Z_1 & -Z_2 \end{array} \right)
\end{array}$$
and $Z_1 = Z_2 = 0$.
Hence, $Z = O$ and ${\rm ad}\ Z = 0$, which leads to
$$InnD(\mathfrak{sl}_2(R)) \cap Der_F(R) = (0).$$
This establishes the result.
QED
\bigskip

We close by making several remarks. The first two point out some results which may be of interest, while the 
last three compare our results and those of \cite{P1}, \cite{P2}.
\bigskip

\noindent
{\bf Remarks}.
\noindent
(1.) Suppose that $R$ is a commutative $F$-algebra. Then it is easy to see that
$\mathfrak{sl}_2(R)$ is finitely generated as a Lie algebra
over $F$ if and only if $R$ is finitely generated as an $F$-algebra.
\bigskip

\noindent
(2.) Suppose that
$F$ is a field of characteristic not $2$ and that
$R$ is an integral domain.
If $(X,H,Y)$ is an $s\ell_2$-triple in $\mathfrak{sl}_2(R)$,
then $FH$ is a MAD. This follows from Theorem 6.
\bigskip

We next give an example in characteristic $0$ where we can use our results to 
conclude the conjugacy of MAD's but where those in \cite{P1}, \cite{P2} cannot be used to draw this conclusion.
\bigskip

\noindent
(3.) Let $F = {\bf Q}$ be the field of rational numbers and put
$$R = F[\xi][ \frac{1}{\xi - \alpha} ]_{\alpha \in F},$$
where $\xi$ is an indeterminate. This is just the localization 
of $F[\xi]$ relative to the monoid generated by the degree $1$ polynomials 
$ \{ \frac{1}{\xi - \alpha} |\alpha \in F\}$. Thus
 $R$ is a principal ideal domain, hence a UFD, so
all of our results are valid for this $R$. 
One knows that $Pic(R)$ is trivial (this is true for any UFD)
and it is easy to see that there is no $F$-rational
point of  $Spec(R)$. Thus our results hold in this case to give the conjugacy of MAD's.  
However the hypothesis of \cite{P1}, \cite{P2} do  not hold in this case so we cannot use the results there to 
obtain the conjugacy of MAD's in this case.
\bigskip

Our next example shows we may use the results of \cite{P1}, \cite{P2} 
in some cases even when we cannot use the ones of this paper. 
\bigskip

\noindent
(4.) Let $F = {\bf Q}$, and we define
$R = {\bf Q} + \xi {\bf R}[\xi]$ as a subring of $S = {\bf R}[\xi]$,
the ring of polynomials in an indeterminate $\xi$ with coefficients
in the field, ${\bf R}$, of real numbers.
Then $R$ is an  integral domain but is not a unique factorization
domain. However the conjugacy theorem for MAD's from \cite{P1}, \cite{P2}
of $\mathfrak{sl}_2(R)$ holds here. Indeed,
 $Pic(R)$ is trivial in this case ( for this see \cite{C}, Remark 3.7). Moreover, there is an $F$-rational point, for example $\xi {\bf R}[\xi]$,
of the scheme $Spec(R)$. Thus the results of \cite{P1}, \cite{P2} apply in this case.
\bigskip

Our final example is where the ring $R$ is not a UFD and also does not staisfy the hypothesis of 
\cite{P1}, \cite{P2}. However we point out how some of our results do hold in this case.
\bigskip

\noindent
(5.) We take $R = F[\xi_1^2,\xi_1^3,\xi_2]$.
Then, $R$ is an integral domain, but not a unique factorization
domain. Furthermore, one knows that $Pic(R)$ is nontrivial (see \cite{C} Example 2.3).
At this moment we cannnot say anything about the conjugacy of MADs.
But we do have a conjugacy of $s\ell_2$-triples under 
the action of the full automorphism group by Theorem 5.

\end{document}